\documentclass[11pt]{article}

\usepackage{amsfonts}
\usepackage{amsmath}
\usepackage{microtype} %omitting bad boxes
\usepackage{multirow}
\usepackage{caption}
\usepackage{subcaption}
\usepackage{graphicx}
\usepackage{color}
\usepackage{epsfig}
\usepackage{graphicx}
\usepackage{subfig}
\usepackage{rotating}
\usepackage{amssymb}
\usepackage{listings}
\usepackage{mathabx} % for widecheck
\newtheorem{lemma}{Lemma}
\newtheorem{theorem}{Theorem}

\newcommand{\Adj}[2]{{\mathrm {Adj}}_{#2}(#1)}
\newcommand{\mrd}{{\mathrm{d}}}
\newcommand{\mre}{{\mathrm{e}}}
\newcommand{\mri}{{\mathrm{i}}}

\newcommand*{\COMMENTS}{}%
\ifdefined\COMMENTS
\newcommand{\comment}[1]{ {\color{red} #1}}% 
\else
\newcommand{\comment}[1]{}
\fi 

\newcommand*{\COMMENTSN}{}%
\ifdefined\COMMENTSN
\newcommand{\commentn}[1]{ {\color{blue} #1}}% 
\else
\newcommand{\commentn}[1]{}
\fi 

\begin{document}

\title{An explicit reconstruction algorithm for the transverse ray transform of a second rank tensor field from three axis data.}
\author{Naeem M. Desai \and William R.B. Lionheart }
\maketitle

\abstract{We give an explicit plane-by-plane filtered back-projection reconstruction algorithm for the transverse ray transform of symmetric second rank tensor fields on Euclidean 3-space, using data from rotation about three orthogonal axes. We show that in the general case two axis data is insufficient but give an explicit reconstruction procedure for the potential case with two axis data.}

\section{Introduction}
\label{Section:Introduction}
The transverse ray transform (TRT) of rank two symmetric tensor fields in three dimensional Euclidean space has recently 
been shown to be of importance in x-ray diffraction strain tomography~\cite{LW}, however currently reconstruction 
algorithms are known only for data from rotations about six axes \cite[Sec 5.1.6]{SBook} or complete data \cite{LW}. 
Given that, in the proposed application, each projection is acquired laboriously using a raster scan, it is advantageous 
to perform the reconstruction using data from a minimum number of axes and making the most of the data collected from 
each projection. In this paper we give an explicit reconstruction formula for three axis data using similar techniques 
to those employed by \cite{LS} for the {\em truncated} transverse ray transform (TTRT). The inversion formula uses 
plane-by-plane filter and back-projection operations familiar from  inversion of the parallel x-ray transform of a scalar 
field.  We go on to show that data from only two rotation axes are insufficient in the general case, but give an explicit 
reconstruction for the potential case with two axis data. We present some numerical results for our reconstruction 
algorithms using simulated data.

\begin{section}{Definitions and notation}
\label{Section:DefinitionsAndNotation}

We denote the complex vector space of symmetric $\mathbb{R}$-bilinear maps $ \mathbb{R}^{3} \times \mathbb{R}^{3} \rightarrow \mathbb{C}$ by
$ S^{2}\mathbb{C}^{3}$ .
The elements of this space are (\textit{complex-valued}) symmetric tensors of second rank on $\mathbb{R}^{3}$. We will identify a complex symmetric tensor $f \in 
S^{2}\mathbb{C}^{3}$  with the $\mathbb{C}$-linear operator, $f : \mathbb{C}^{n} \rightarrow \mathbb{C}^{n},$ as usual  $\langle f\xi,\eta\rangle = \langle f\eta,\xi\rangle = f(\xi,\eta)$ where $\eta,\xi \in \mathbb{R}^{3}$. 

Defining the Schwartz space of rapidly decreasing functions  $\mathcal{S}(\mathbb{R}^n)$ on $n$-dimensional Euclidean ($n=2$ or $3$) space in the usual way, we extend this definition to (complex) vector fields $\mathcal{S}(\mathbb{R}^n;\mathbb{C}^n)$ and complex symmetric tensor fields as $\mathcal{S}(\mathbb{R}^n;S^2\mathbb{C}^n)$. The choice of  Schwartz spaces and the use of complex vectors and tensors is convenient as we will rely heavily on the Fourier transform. Extension to Sobolev spaces for compactly supported fields follows using the usual apparatus as applied to scalar ray transforms \cite{Natterer}.
We  define the Fourier transform 
$F : \mathcal{S}(\mathbb{R}^{3}) \rightarrow \mathcal{S}(\mathbb{R}^{3})$, by
$$ F\left[f\right] = \hat{f}(y) = (2\pi)^{-3/2}\int_{\mathbb{R}^{3}} \mre^{-\mri\langle y,x\rangle}f(x) \,\mrd x.$$

Given an orthonormal basis $\left(e_1,e_2,e_3\right)$ of $\mathbb{R}^{3}$, a tensor $f \in S^{2}\mathbb{C}^{3}$ can be 
represented by the symmetric $3 \times 3$ matrix $(f_{jk}), f_{jk} = f(e_j,e_k).$ The Hermitian scalar product on 
$S^{2}\mathbb{C}^{3}$  can be written as $\langle f,g\rangle = \sum^{3}_{j,k=1} f_{jk}\bar{g}_{jk}$ independently 
of the choice of an orthonormal basis. Since only orthonormal bases will be used in this paper, we do not distinguish between 
co- and contravariant tensors.

We will need the partial Fourier transform $F_V : \mathcal{S}(\mathbb{R}^{3}) \rightarrow \mathcal{S}(\mathbb{R}^{3})$, 
for any $k$-dimensional vector subspace $V \subset \mathbb{R}^{3}$. Given Cartesian coordinates $(x_1,x_2,x_3)$ in $\mathbb{R}^{3}$ such that
\\$V = \{x|x_{k+1} = .... = x_3 = 0\}.$  the partial Fourier Transform can be written as
$$ F_V\left[f\right] = \hat{f}(y_1,...,y_k,x_{k+1},...,x_3) = (2\pi)^{-k/2}\int_{\mathbb{R}^{k}}\mre^{-\mri(y_1x_1+...+y_kx_k)}f(x)\,\mrd x_1...\,\mrd x_k. $$
This result is independent of the choice of orthonormal coordinates and satisfies the commutation law $F = F_VF_{V^{\bot}} = F_{V^{\bot}}F_V.$

The {\em oriented} lines in $\mathbb{R}^{3}$ can be parameterized  by points of the tangent bundle of the unit sphere $\mathbb{S}^2$ in $\mathbb{R}^3$
$$ T\mathbb{S}^{2} = \{(\xi,x) \in \mathbb{R}^{3} \times \mathbb{R}^{3} \mid |\xi| = 1, \langle\xi,x\rangle = 0\} 
\subset \mathbb{R}^{3} \times \mathbb{R}^{3}.$$
An oriented line $l \subset \mathbb{R}^{3}$ is uniquely represented as $l = \{x + t\xi \mid t \in \mathbb{R}\}$ 
with $(\xi,x) \in T\mathbb{S}^{2}$. The Schwartz space $\mathcal{S}(T\mathbb{S}^{2})$ is defined as in \cite{LS}.

Our main object of interest is the  Transverse Ray Transform (TRT) of symmetric rank two tensor fields
$$J:\mathcal{S}(\mathbb{R}^{3};S^{2}\mathbb{C}^{3}) \rightarrow \mathcal{S}(T\mathbb{S}^{2};S^{2}\mathbb{C}^{3}),$$ which is
defined by
\begin{equation}
Jf(\xi,x) = \int_{-\infty}^{\infty} P_{\xi}f(x + t\xi) \,\mrd t,
\label{eqn:TRT}
\end{equation}
where $P_{\xi}: S^{2}\mathbb{C}^{3} \rightarrow S^{2}\mathbb{C}^{3}$ is the orthogonal projection onto the subspace
\\$\{f \in S^{2}\mathbb{C}^3 \mid f\xi = 0 \}.$ 

For example, for an orthonormal basis of the form $(e_1,e_2,e_3 = \xi)$, the projection 
is expressed by
\begin{equation}
P_{\xi}f = \left(\begin{array}{ccc} f_{11} & f_{12} & 0 \\ f_{12} & f_{22} & 0 \\ 0 & 0 & 0 \end{array} \right).
\label{eqn:J_project}
\end{equation}

By  contrast the scalar x-ray transform (or simply the ray transform) of a function $f\in \mathcal{S}(\mathbb{R}^3)$ is defined by 
$$ Xf(\xi,x) = \int_{-\infty}^{\infty} f(x + t\xi) \,\mrd t.$$
We also have the Longitudinal Ray Transform (LRT), 
$$I:\mathcal{S}(\mathbb{R}^{3};\mathbb{C}^{3}) \rightarrow \mathcal{S}(T\mathbb{S}^{2}), \quad
I:\mathcal{S}(\mathbb{R}^{3};S^{2}\mathbb{C}^{3}) \rightarrow \mathcal{S}(T\mathbb{S}^{2}),$$  defined on vector
and tensor fields respectively by
\begin{equation}
 If(\xi,x) = \int_{-\infty}^{\infty} \langle f(x + t\xi),\xi \rangle \,\mrd t, \quad 
 If(\xi,x) = \int_{-\infty}^{\infty} \langle f(x + t\xi)\xi,\xi \rangle \,\mrd t.
\label{eqn:LT}
\end{equation}
%Plane not hyperplane!
The x-ray transform and LRT can also be defined on a plane in $\mathbb{R}^{3}$. For $\eta \in \mathbb{S}^{2}$, let 
$\eta^{\bot} = \{\xi \in \mathbb{R}^{3} \mid \langle \xi, \eta \rangle = 0 \}$,  $\mathbb{R}\eta =  \{s\eta \mid s \in \mathbb{R} \}$, $\eta^{\bot}_{\mathbb{C}}$ be the 
complexification of $\eta^{\bot}$ and \newline $\mathbb{S}_{\eta}^{1} = \{\xi \in \eta^{\bot}\mid |\xi| = 1\} $ 
be the unit sphere in $\xi^{\bot}$.
Given $s \in \mathbb{R}$, let $s\eta + \eta^{\bot}$ be the plane through $s\eta$ parallel to $\eta^{\bot}$ and
$\iota_{s,\eta} : s\eta + \eta^{\bot} \subset \mathbb{R}^{3}$ be the identical embedding. The family of oriented lines in the 
plane $s\eta + \eta^{\bot} $ is parameterized by points of the manifold \newline
$T\mathbb{S}_{\eta}^{1} = \{(\xi,x) \mid \xi \in \mathbb{S}^{1}_{\eta}, x \in \eta^{\bot}, \langle \xi,x \rangle = 0 \}$
such that a point $(\xi,x) \in T\mathbb{S}_{\eta}^{1}$ corresponds to the line $\{s\eta + x + t\xi \mid t \in \mathbb{R} \}.$
We define the X-ray transform on a plane $X_{\eta,s}f(\xi,x)$ and the LRT on the plane $s\eta + \eta^{\bot}$
\begin{eqnarray}
X_{\eta,s}:& \mathcal{S}(s\eta + \eta^{\bot};\mathbb{C}) \rightarrow \mathcal{S}(T\mathbb{S}^{1}_{\eta}) \\
I_{\eta,s}:& \mathcal{S}(s\eta + \eta^{\bot};\eta^{\bot}_{\mathbb{C}}) \rightarrow \mathcal{S}(T\mathbb{S}^{1}_{\eta})\\
 I_{\eta,s}:& \mathcal{S}(s\eta + \eta^{\bot};S^2\eta^{\bot}_{\mathbb{C}}) 
\rightarrow \mathcal{S}(T\mathbb{S}^{1}_{\eta})
\end{eqnarray}
by the following formulae
\begin{eqnarray}
X_{\eta,s}f(\xi,x) &= \int_{-\infty}^{\infty}  f(s\eta+x+t\xi) \, \,\mrd t, \\
I_{\eta,s}f(\xi,x) &= \int_{-\infty}^{\infty} \langle f(s\eta+x+t\xi),\xi \rangle \, \,\mrd t, \\
I_{\eta,s}f(\xi,x) &= \int_{-\infty}^{\infty} \langle f(s\eta+x+t\xi)\xi,\xi \rangle \ \,\mrd t
\label{eqn:LT_on_plane}
\end{eqnarray}
respectively. One can see that operators (\ref{eqn:LT}) and (\ref{eqn:LT_on_plane}) are related. 
If $f \in \mathcal{S}(\mathbb{R}^{3};S^2\mathbb{C}^{3})$ and $\iota^{*}_{\eta,s}f$ is the slice of $f$ 
by the plane $s\eta + \eta^{\bot}$, then \\$I_{\eta,s}(\iota^{*}_{\eta,s}f)(\xi,x) = If(\xi,s\eta + x), \text{for} \ 
(\xi,x) \in T\mathbb{S}^{1}_{\eta}$.

A typical experimental situation would involve rotation of the specimen (or equivalently the source and detector) about some finite collection of axes. 
In the scalar case of the  x-ray transform rotation about one axis is sufficient as the problem reduces to the Radon transform in the plane, for a family of planes normal to the rotation  axis. Consider first the case $n=2$ in which $X$ is identical to the Radon transform. In this case the formal adjoint $B=X^*: \mathcal{S}(T\mathbb{S}^{1}) \rightarrow C^{\infty}(\mathbb{R}^{2})$, 
the back-projection operator, is well defined and given by
\begin{equation}
B\phi(x) =  \frac{1}{2 \pi}\int_{\mathbb{S}^1} \phi(\xi,\langle\xi, x\rangle) d \xi
\label{eqn:BackProjection}
\end{equation}
for $\phi \in \mathcal{S}(T\mathbb{S}^{1})$.
We then have the inversion formula \cite{Natterer} for data $\phi(\xi,x)=X_{\eta,s}f(\xi,x)$ in the range of $X$
\begin{equation}
\label{eqn:FBP}
F[f(x)](y) = | y | F[ B\phi(x)],
\end{equation}
which means that inversion is performed by a  {\em ramp filter} (also known as a Riesz potential)  applied to the back-projected data. 

This operation can be performed {\em slice by slice} to invert the x-ray transform for $n=3$, in which case data is needed only for $\xi\in \eta^\perp$ for some fixed {\em rotation axis} $\eta \in \mathbb{S}^2$. In this case the  slice-by-slice back-projection  operator 
$B_{\eta}: \mathcal{S}(\mathbb{R} \times T\mathbb{S}_{\eta}^{1}) \rightarrow C^{\infty}(\mathbb{R}^{3})$.
% I think it is obvious
% defined by:
%\begin{equation}
%(B_{\eta}\phi)(x) = 
%\frac{1}{2 \pi}\int_{S^1} \phi(\xi(\theta),\langle\xi(\theta), x\rangle),\mrd \theta, \quad \text{for} \quad x \in \eta^{\bot},
%\label{eqn:BackProjection_eta}
%\end{equation}
%where $\xi(\theta) =  \cos\theta e_1 + \sin\theta e_2$ and $\xi^{\bot}(\theta) =  \cos\theta e_2 -\sin\theta e_1$ with respect to 
%an orthonormal basis $(e_1,e_2)$ of $\eta^{\bot}.$ The value of $B_{\eta}\phi$ at a point $x \in \mathbb{R}^{3}$ is just the 
%average of the function $\phi$ over all lines passing through $x$ and orthogonal to $\eta$.
 The reconstruction formula (\ref{eqn:FBP}) becomes
\begin{equation}
\label{eqn:slicebyslice}
F_{\eta^\perp}[f(x)](s,y) = | \Pi_\eta y | F_{\eta^\perp}[ B_\eta\phi(x)],
\end{equation}
where $s= \langle x, \eta \rangle$.

Notice that the component $\langle \eta, Jf(\xi,x)\eta\rangle = X\langle \eta, f\eta\rangle(x,\xi)$, for $\xi \in \eta^\perp$ is simply the scalar x-ray transform in the plane through $x$ normal to $\eta$. As observed in \cite[Sec 5.1.6]{SBook} this component can be reconstructed using any inversion formula for the planar Radon transform inversion plane by plane, including the one given in (\ref{eqn:slicebyslice}). Choosing six rotation axes $\eta_i$ so that the outer products $\eta \eta^T$ are linearly independent in $S^2\mathbb{R}^3$ recovers $f$ everywhere. As mentioned in the introduction this procedure is likely to be time-consuming as rotations are performed about six axes and yet for each ray only one measurement (out of a possible three) is used. The aim of this paper is to show, via a constructive inversion procedure, that $f$ can be determined uniquely from the data $Jf(x,\xi)$ for $\xi \in e_i^\perp$, $i=1,...,3$. Of course the diagonal elements $f_{ii}$ are already determined as above so our main task is to provide a reconstruction procedure for the off-diagonal elements.

For a given rotation axis $\eta$ and direction $\xi\in \eta^\perp$ we have in addition to the `axial' component 
$\langle \eta, Jf(\xi,x)\eta\rangle$ the `non axial' components $\langle \zeta, Jf(\xi,x)\eta\rangle$ where $\zeta = \xi\times\eta$ and $\langle \zeta, Jf(\xi,x)\zeta\rangle$ with which to reconstruct the off-diagonal (in a basis including $\eta$) elements of $f$.

\end{section}

\begin{section}{Relations between transforms}

The aim of this section is to write the non-axial components of $Jf(\xi,x)$, $\xi \in \eta^\perp$ in terms of longitudinal ray transforms on transaxial planes.

We have  $P_{\xi}(f) = \Pi_{\xi}f\Pi_{\xi}$ where $\Pi_{\xi}$ is the orthogonal
projection matrix onto $\xi^{\bot},$ and the `off diagonal' component can be expressed as
$$\langle(Jf)(\xi,x)\eta,\xi\times\eta \rangle = 
\int_{-\infty}^{\infty} \langle \Pi_{\xi}f(x + t\xi)\Pi_{\xi}\eta,\xi\times\eta \rangle \,\mrd t $$
$$ = \int_{-\infty}^{\infty} \langle f(x + t\xi)\eta,\xi\times\eta \rangle \,\mrd t = 
\int_{-\infty}^{\infty} \langle \eta \times f(x + t\xi)\eta,\xi \rangle \,\mrd t. $$

Since the vector field $\eta \times f\eta$ is orthogonal to $\eta$, its restriction to every plane $s\eta + \eta^{\bot}$ can be
considered as a vector field on the plane, i.e., 
$(\eta\times f\eta)|_{s\eta + \eta^{\bot}} \in \mathcal{S}(s\eta + \eta^{\bot};\eta^{\bot}_{\mathbb{C}}).$ As this is then contracted with the ray direction $\xi$ we have
\begin{equation}
I_{\eta,s}((\eta\times f\eta)|_{s\eta + \eta^{\bot}})(\xi,x) = \langle(Jf)(\xi,x)\eta,\xi\times\eta \rangle \quad \text{for}
\quad (\xi,x) \in T\mathbb{S}^{1}_{\eta}.
\label{eqn:offDiag3}
\end{equation}
As we have seen in (\ref{eqn:J_project}), that the TRT depends only on the projection normal to the direction of the ray. Now working in $(\eta,\zeta = \xi\times\eta,\xi)$ coordinates, let us consider, 
$\langle Jf(\xi,x)\zeta,\zeta \rangle$, which can be transformed into
\begin{equation}
\langle Jf(\xi,x)\zeta,\zeta \rangle = \int_{-\infty}^{\infty} \langle f(x + t\xi)\zeta,\zeta \rangle \ \,\mrd t.
\label{eqn:Jf_zeta}
\end{equation}
Let us parameterize $\xi$ in the usual sense as
$$\xi = \left( \begin{array}{c} \cos\theta \\ \sin\theta \\ 0 \end{array} \right), \quad \text{so,} \quad 
\zeta = \left( \begin{array}{c} -\sin\theta \\ \cos\theta \\ 0 \end{array} \right). $$
Since $\zeta \in \xi^{\bot}$ we can calculate
$$ \langle \zeta, f\zeta\rangle = (-\sin\theta,\cos\theta,0)
\left(\begin{array}{ccc} f_{11} & f_{12} & 0 \\ f_{12} & f_{22} & 0 \\ 0 & 0 & 0 \end{array} \right) 
\left( \begin{array}{c} -\sin\theta \\ \cos\theta \\ 0 \end{array} \right) $$ $$ = f_{22}\cos^2\theta - 
2f_{12}\cos\theta\sin\theta + f_{11}\sin^2\theta $$ $$ = (\cos\theta,\sin\theta,0)
\left(\begin{array}{ccc} f_{22} & -f_{12} & 0 \\ -f_{12} & f_{11} & 0 \\ 0 & 0 & 0 \end{array} \right) 
\left( \begin{array}{c} \cos\theta \\ \sin\theta \\ 0 \end{array} \right) $$ $$ = \xi \Adj{f}{e_3^\perp} \xi,$$ 
where $\Adj{f}{e_3^\perp}$ denotes the adjugate matrix  of the slice of $f$ restricted to the plane; of course this is nothing other than the  conjugation of $P_{e_3}f$ with a right angle rotation about the $e_3$ axis.
Hence using the above, we see that for $\xi\in\eta^\perp$
\begin{equation}
\langle Jf(\xi,x)\zeta,\zeta \rangle = \int_{-\infty}^{\infty} \langle f(x + t\xi)\zeta,\zeta \rangle \ \,\mrd t $$ $$
 = \int_{-\infty}^{\infty} \langle \Adj{f}{\eta^\perp}(x+ t\xi)\xi,\xi \rangle = \langle I\Adj{f}{\eta^\perp}(\xi,x)\xi,\xi \rangle,
\label{eqn:Jf=IAdjf}
\end{equation}
which is the LRT transform of the two dimensional adjugate of the slice of $f$ restricted to the plane. 
We notice that this is also exactly the transverse ray transform in the planar case.
The results of this section can be summarized in the following Lemma:

\begin{lemma}\label{Lemma3.1}
Let $f \in \mathcal{S}(\mathbb{R}^{3};S^2\mathbb{C}^{3})$ be a symmetric tensor field. The equations
\begin{equation}
I_{\eta,s}((\eta\times f\eta)|_{s\eta + \eta^{\bot}}) = (J^{1}_{\eta,s}f),
\label{eqn:Lemma3.1a}
\end{equation}
\begin{equation}
I_{\eta,s}\left(\Adj{\iota^{*}_{\eta,s}f}{\eta^\perp}\right) = (J^{2}_{\eta,s}f),
\label{eqn:Lemma3.1b}
\end{equation}
hold for every $s \in \mathbb{R}$ and $\eta \in \mathbb{S}^{2}$, where
$$ (J^{1}_{\eta,s}f) = \langle(Jf)(\xi,s\eta + x)\eta,\xi\times\eta \rangle, 
\quad (J^{2}_{\eta,s}f) = \langle(Jf)(\xi,s\eta + x)\zeta,\zeta\rangle.$$
\label{lemma:section3}
\end{lemma}

 In the next section, we will transform (\ref{eqn:Lemma3.1a}) and (\ref{eqn:Lemma3.1b})
to algebraic equations by applying the Fourier transform to back-projected data.

\end{section}

\begin{section}{Main algebraic equations}

\begin{subsection}{Curl Components of Tensor and Vector Fields}

We can now transform (\ref{eqn:Lemma3.1a}) and (\ref{eqn:Lemma3.1b}) to algebraic equations by applying the Fourier transform.
We only require what \cite{SBook} refers to as  tangential component $\tau g \in C^{\infty}(\mathbb{R}^{2})$ of a vector field $g \in C^{\infty}(\mathbb{R}^{2};\mathbb{C}^{2})$,
which is defined by
\begin{equation}
(\tau g)(y) = \langle g(y),y^{\bot} \rangle.
\label{eqn:TangComp}
\end{equation}
Here the vector $y^{\bot}$ is the result of rotating $y$ by $\pi/2$ in the positive direction. Of course, one can understand 
(\ref{eqn:TangComp}) as the $2D$ curl of a vector field in Fourier (frequency) space. 
The manifold $T\mathbb{S}^{1}$
can be identified with $\mathbb{R}\times\mathbb{S}^{1}$ by the diffeomorphism $(p,\xi) \mapsto (\xi,p\xi^{\bot})$ for
$(p,\xi) \in \mathbb{R}\times\mathbb{S}^{1}$. Therefore the derivative 
$\frac{\partial}{\partial p} : \mathcal{S}(T\mathbb{S}^{1}) \rightarrow \mathcal{S}(T\mathbb{S}^{1})$ is well defined.
For a vector field $f \in \mathcal{S}(\mathbb{R}^{2};\mathbb{C}^{2})$, the tangential component of the Fourier Transform $F[f]$
is recovered by the LRT, $If$, by the formula
\begin{equation}
\tau F[f] = \frac{\mri}{2}|y|F\left[B\left(\frac{\partial(If)}{\partial p}\right)\right].
\label{eqn:FourierTangComp}
\end{equation}

We see in  \cite{LS} and \cite{SliceS}, the tangential component, $\tau g \in C^{\infty}(\mathbb{R}^2)$, of a tensor 
field $g \in C^{\infty}(\mathbb{R}^2;S^2\mathbb{C}^2)$ is defined by
\begin{equation}
(\tau g)(y) = |y|^2\mathrm{tr\,}g - \langle g(y)y,y \rangle.
\label{eqn:TangCompTens}
\end{equation}
This is exactly the Fourier transform of the single unique non-zero component of the 
compatibility tensor of Barr\'{e} de Saint-Venant in the plane 
\begin{equation}
  W(f)=  \frac{\partial^2 f_{11}}{\partial x_2^2}
   - 2\frac{\partial^2 f_{12}}{\partial x_1 \partial x_2}
   + \frac{\partial^2 f_{22}}{\partial x_1^2},	
\end{equation}
which is also sometimes described as the curl of a symmetric tensor field.
For $f \in \mathcal{S}(\mathbb{R}^2;S^2\mathbb{C}^2)$, the tangential component of the Fourier transform $F[f]$ is 
recovered from the LRT, $If,$ as
\begin{equation}
\tau F[f] = \frac{1}{2}|y|^3F[B(If)].
\label{eqn:FourierTangCompTens}
\end{equation}
For $\phi \in \mathcal{S}(T\mathbb{S}^1)$, the function $B\phi(x)$ is $C^{\infty}$-smooth and bounded on 
$\mathbb{R}^2$ but does not decay fast enough to be in the Schwartz class. Thus we understand the Fourier 
transform in the distribution sense in (\ref{eqn:FourierTangComp}) and (\ref{eqn:FourierTangCompTens}).

\end{subsection}

\begin{subsection}{Derivation of System of Equations}

Let $f \in \mathcal{S}(\mathbb{R}^{3};S^{2}\mathbb{C}^{3})$ be a symmetric tensor field and denote by 
$f^{\prime} = F_{\eta^{\bot}}[f] \in \mathcal{S}(\mathbb{R}^{3};S^{2}\mathbb{C}^{3})$ the partial Fourier transform of $f$.
For any $s \in \mathbb{R}$, the restriction of the vector field $\eta\times f^{\prime}\eta$ to the plane $s\eta + \eta^{\bot}$ 
coincides with the 2D Fourier transform of $(\eta\times f\eta)|_{s\eta + \eta^{\bot}}$. This is 
$(\eta\times f^{\prime}\eta)|_{s\eta + \eta^{\bot}} = F_{\eta^\perp}[(\eta\times f\eta)|_{s\eta + \eta^{\bot}}].$ \\ \\
We then apply formula (\ref{eqn:FourierTangComp}) to the vector field $(\eta\times f\eta)|_{s\eta + \eta^{\bot}}$, giving
\begin{equation}
\tau((\eta\times f\eta)|_{s\eta + \eta^{\bot}})(s\eta + y) = \frac{\mri}{2}|y|F_{\eta^{\bot}}\left[B_{\eta}\left(
\frac{\partial(I_{\eta,s}((\eta\times f\eta)|_{s\eta + \eta^{\bot}}))}{\partial p}\right)\right]\quad \text{for}\quad y\in \eta^{\bot},
\label{eqn:Lambda1}
\end{equation}
Using (\ref{eqn:Lemma3.1a}), we can transform (\ref{eqn:Lambda1}) giving
\begin{equation}
\tau((\eta\times f^{\prime}\eta)|_{s\eta + \eta^{\bot}})(s\eta + y) = \frac{\mri}{2}|y|F_{\eta^\perp}\left[\left(B_{\eta}
\frac{\partial(J^{1}_{\eta}f)}{\partial p}\right)(s\eta + x)\right],
\label{eqn:Lambda2}
\end{equation}
%where $F_{x\rightarrow y}$ is the Fourier transform on the plane $\eta^{\bot}$.
Note that (\ref{eqn:TangComp}) gives
\begin{equation}
\tau((\eta\times f^{\prime}\eta)|_{s\eta + \eta^{\bot}})(s\eta + y) = \langle\eta\times f^{\prime}(s\eta + y)\eta,\eta\times y
\rangle = \langle f^{\prime}(s\eta + y)\eta,y\rangle.
\label{eqn:TangCompEtaCrossfEta}
\end{equation}
Upon substitution of (\ref{eqn:TangCompEtaCrossfEta}) into the LHS of (\ref{eqn:Lambda2}) and applying the 
one-dimensional Fourier transform $F_{\mathbb{R}\eta}$  taking $s$ to $\sigma$ gives
\begin{equation}
\langle \hat{f}(\sigma\eta + y^{\prime})\eta,y^{\prime}\rangle = \frac{\mri}{2}|y^{\prime}|
%F_{s\eta + x^{\prime}\rightarrow \sigma\eta+y^{\prime}}
F
\left[\left(B_{\eta}\frac{\partial(J^{1}_{\eta}f)}{\partial p}\right)(s\eta + x^{\prime})\right] \quad \text{for} \quad 
y^{\prime} \in \eta^{\bot},
\label{eqn:OneDimenFourierEtaxfEta}
\end{equation}
where $\hat{f}$ is the three-dimensional Fourier transform  $F[f]$. Since $y^{\prime} \in \eta^{\bot}$ and $\sigma \in \mathbb{R}$, we let
$ y = \sigma\eta + y^{\prime}$, where $y^{\prime} = \Pi_{\eta}y$. Hence the previous
formula (\ref{eqn:OneDimenFourierEtaxfEta}) can be written as
\begin{equation}
\langle \hat{f}(y)\eta,\Pi_{\eta}y\rangle = \frac{\mri}{2}|\Pi_{\eta}y|F
\left[\left(B_{\eta}\frac{\partial(J^{1}_{\eta}f)}{\partial p}\right)(x)\right],
\label{eqn:Lambda3}
\end{equation}
%where $F_{x\rightarrow y}$ is the three-dimensional Fourier transform.
Note that (\ref{eqn:Lambda3}) is identical to the 
off-diagonals for the TTRT operator case in \cite{LS}. Moreover this will just reconstruct the solenoidal part
of the off-diagonals since the Fourier transform interweaves with the solenoidal part. 
%Again we consider $\eta \in \mathbb{S}^2$ to be a fixed vector and let $f^{\prime} = F_{\eta^{\bot}}[f]$ be the 
%partial Fourier transform of a symmetric tensor field $f \in \mathcal{S}(\mathbb{R}^3;S^2\mathbb{C}^3)$.
 For any 
$s \in \mathbb{R},$ the slice $ \iota^{*}_{\eta,s}f^{\prime}$ coincides with the $2D$ Fourier transform of the slice 
$\iota^{*}_{\eta,s}f$, i.e., $\iota^{*}_{\eta,s}f^{\prime} = F_{\eta^\perp}[\iota^{*}_{\eta,s}f]$, where the Fourier 
transform on the plane $s\eta + \eta^{\bot}.$ 
%With the aid of Lemma \ref{Lemma3.1} we relabel  
%$ h =  \Adj{\iota^{*}_{\eta,s}f}{\eta^\perp}$ and with similar reasoning $\overline{\iota^{*}_{\eta,s}f^{\prime}} 
%= F[\overline{\iota^{*}_{\eta,s}f}]$, since we are considering the adjugate of the slice of $f$ restricted to the 
%plane $s\eta + \eta^{\bot}.$ 
Henceforth, we refer to the adjugate of the slice of $f^{\prime}$ restricted to the plane as
$\Adj{\iota^{*}_{\eta,s}f^{\prime}}{\eta^\perp} = h^{\prime}.$ 
Upon application of formula (\ref{eqn:FourierTangCompTens}) to $h^{\prime},$ we see
%application of formula (\ref{eqn:TangCompTens}) to the 
%tensor field $\iota^{*}_{\eta,s}f \in \mathcal{S}(s\eta + \eta^{\bot};S^2\eta^{\bot}_{\mathbb{C}})$, we now see
\begin{equation}
[\tau(h^{\prime})](s\eta + y) = \frac{1}{2}|y|^3F_{\eta^{\bot}}[B_{\eta}(I_{\eta,s}(\Adj{\iota^{*}_{\eta,s}f}{\eta^\perp}))] \ \ 
\text{for} \ \ y \in \eta^{\bot}.
\label{eqn:Mu1}
\end{equation}
Using Lemma \ref{lemma:section3} we can rewrite the above as
\begin{equation}
[\tau(h^{\prime})](s\eta + y) = \frac{1}{2}|y|^3F_{\eta^\perp}[B_{\eta}(J_{\eta}^2f)(s\eta + x)] \ \ 
\text{for} \ \ y \in \eta^{\bot},
\label{eqn:Mu2}
\end{equation}
 Now, we apply formula 
(\ref{eqn:TangCompTens}) to the field $g = h^{\prime} \in \mathcal{S}(s\eta + 
\eta^{\bot};S^2\eta^{\bot}_{\mathbb{C}})$ to give
\begin{equation}
[\tau(h^{\prime})](s\eta + y) = |y|^2\mathrm{tr\,}h^{\prime}(s\eta +y) - \langle h^{\prime}(s\eta + y)y,y \rangle \ \ 
\text{for} \ \ y \in \eta^{\bot}.
\label{eqn:Mu3}
\end{equation}
Next, substitute (\ref{eqn:Mu3}) into (\ref{eqn:Mu2}) to give
\begin{equation}
|y|^2\mathrm{tr\,}h^{\prime}(s\eta + y) - \langle h^{\prime}(s\eta + y)y,y \rangle = 
\frac{1}{2}|y|^3F_{\eta^\perp}[B_{\eta}(J_{\eta}^2f)(s\eta + x)]\ \ \text{for} \ \ y \in \eta^{\bot}.
\label{eqn:Mu4}
\end{equation}
By applying the one-dimensional Fourier transform on $\mathbb{R}\eta$ to the above, we obtain
\begin{equation}
|y^{\prime}|^2\mathrm{tr\,}\hat{h}(\sigma\eta + y^{\prime}) - \langle \hat{h}(\sigma\eta + y^{\prime})y^{\prime},y^{\prime}
\rangle = \frac{1}{2}|y^{\prime}|^3F
%_{s\eta + x^{\prime} \rightarrow \sigma\eta + y^{\prime}}
[B_{\eta}(J_{\eta}^2f)(s\eta + x^{\prime})],
\label{eqn:Mu5}
\end{equation}
for $y \in \eta^{\bot}.$ As before, employing a change of variables, $y = \sigma\eta + y^{\prime},$ transforms the 
above to
\begin{equation}
|\Pi_{\eta}y|^2\mathrm{tr\,}\hat{h}(y) - \langle \hat{h}(y)\Pi_{\eta}y,\Pi_{\eta}y\rangle = 
\frac{1}{2}|\Pi_{\eta}y|^3F [B_{\eta}(J_{\eta}^2f)(x)], \ \ \text{for} \ \ y \in \mathbb{R}^3.
\label{eqn:Mu6}
\end{equation}

In the following statement, the results are summarized.

\begin{lemma}
Let $\widehat{f}$ be the $3D$ Fourier transform of a symmetric tensor field $f 
\in \mathcal{S}(\mathbb{R}^{3};S^{2}\mathbb{C}^{3})$. For a unit vector $\eta \in \mathbb{S}^{2}$, the following 
equations hold with the additional condition that $\widehat{h} \in \mathcal{S}(\mathbb{R}^{2};S^{2}\mathbb{C}^{2}),$ is 
defined to be the $2D$ adjugate of $f$ restricted to the plane
\begin{equation}
\langle \widehat{f}(y)\eta,\Pi_{\eta}y\rangle = \lambda_{\eta}(y)\quad \text{and}
\label{eqn:Lemma4.1a}
\end{equation}
\begin{equation}
|\pi_{\eta}y|^2\mathrm{tr\,}\hat{h}(y) - \langle \hat{h}(y)\Pi_{\eta}y,\Pi_{\eta}y\rangle = \mu_{\eta}(y),
\label{eqn:Lemma4.2a}
\end{equation}
hold on $\mathbb{R}^{3}$, with the right hand sides defined by
\begin{equation}
\lambda_{\eta}(y) = \frac{\mri}{2}|\Pi_{\eta}y|F
\left[\left(B_{\eta}\frac{\partial(J^{1}_{\eta}f)}{\partial p}\right)(x)\right],
\label{eqn:Lemma4.1b}
\end{equation}
\begin{equation}
\mu_{\eta}(y) = \frac{1}{2}|\Pi_{\eta}y|^3F [B_{\eta}(J_{\eta}^2f)(x)].
\label{eqn:Lemma4.2b}
\end{equation}
\label{lemma:section4}
\end{lemma}
 The partial derivative 
$\frac{\partial}{\partial p}: \mathcal{S}(\mathbb{R}\times T\mathbb{S}^{1}_{\eta}) \rightarrow
\mathcal{S}(\mathbb{R}\times T\mathbb{S}^{1}_{\eta})$ is defined with the help of the diffeomorphism 
$\mathbb{R}^{2} \times \mathbb{S}^{1}_{\eta} \rightarrow \mathbb{R} \times \mathbb{R} \times T\mathbb{S}^{1}_{\eta},
(s,p,\xi) \mapsto (s,\xi,p\xi\times\eta)$.
Given the data $Jf|_{\eta^\perp},$ right-hand sides $\lambda_{\eta}(y)$ and $\mu_{\eta}(y)$ of equations (\ref{eqn:Lemma4.1a}) 
- (\ref{eqn:Lemma4.2a}) can be effectively recovered by formulas (\ref{eqn:Lemma4.1b}) - (\ref{eqn:Lemma4.2b}).

Consider the case where $\eta = \eta_{1} = (1,0,0).$ So $\Pi_{\eta}y = (0,y_2,y_3), 
\widehat{f}\eta = (\widehat{f}_{11},\widehat{f}_{12},\widehat{f}_{13}).$ To abbreviate formulas further, let us denote $\lambda_{\eta_{i}}$ by
$\lambda_{i}$ and $\mu_{\eta_{i}}$ as $\mu_{i}.$ We use the orthonormal basis vectors $e_1,e_2$ and $e_3$ for $\eta$ 
(i.e. $\eta_{1} = e_{1}, \eta_{2} = e_{2}$ and $\eta_{3} = e_{3}$). With the aid of (\ref{eqn:Lemma4.1a}) and 
(\ref{eqn:Lemma4.1b}), we obtain a system of equations
\begin{equation}
\left(\begin{array}{ccc} y_2 & y_3 & 0 \\ y_1 & 0 & y_3 \\ 0 & y_1 & y_2 \end{array} \right)
\left(\begin{array}{c} \widehat{f}_{12}\\ \widehat{f}_{13}\\ \widehat{f}_{23} \end{array} \right) =  
\left(\begin{array}{c} \lambda_{1}\\ \lambda_{2} \\ \lambda_{3} \end{array} \right).  
\label{eqn:SysLambda}
\end{equation}
The system of equations (\ref{eqn:SysLambda}) can be solved to give
\begin{eqnarray}
\widehat{f}_{12} &=& \frac{\lambda_1}{2y_2}+\frac{\lambda_2}{2y_1}-\frac{\lambda_3y_3}{2y_1y_2}, \\
\widehat{f}_{13} &=& \frac{\lambda_1}{2y_3}+\frac{\lambda_3}{2y_1}-\frac{\lambda_2y_2}{2y_1y_3},\\
\widehat{f}_{23} &=& \frac{\lambda_2}{2y_3}+\frac{\lambda_3}{2y_2}-\frac{\lambda_1y_1}{2y_2y_3}.
\label{eqn:SolnLambda}
\end{eqnarray}
The result is summarized in the following theorem

\begin{theorem}
A symmetric tensor field $f \in \mathcal{S}(\mathbb{R}^{3};S^2\mathbb{C}^{3})$ is uniquely determined by the data $Jf(\xi,x)$ for $\xi \in \eta_i^\perp$, $i=1,2,3$ where $(\eta_1,\eta_2,\eta_3)$ forms an orthogonal basis.
\end{theorem}

\end{subsection}
\end{section}

\section{Alternative formulae for diagonal components}
\label{section:alternate}

While the diagonal components $f_{ii}$ are easily determined as we have seen, there is an alternative more complicated procedure to recover them. As this uses different data it can also be viewed as a compatibility condition on the three axis data.

Consider (\ref{eqn:Lemma4.2a}) and (\ref{eqn:Lemma4.2b}). When $\eta = e_1$, we have tr$\hat{h} =  
\hat{f}_{22} + \hat{f}_{33},$ and
$$ \hat{h} = \left(\begin{array}{ccc} 0 & 0 & 0 \\ 0 & \hat{f}_{33} & -\hat{f}_{23} \\ 0 & -\hat{f}_{23} & \hat{f}_{22} 
\end{array} \right). $$
In the same manner as above ($\lambda_i$), we achieve a system of equations for $\mu_{i}$
\begin{equation}
(y_2^2 + y_3^2)(\hat{f}_{22} + \hat{f}_{33}) - (y_2^2\hat{f}_{33} - 2y_2y_3\hat{f}_{23} + 
y_3^2\hat{f}_{22}) = \mu_{1}, $$ $$ 
(y_1^2 + y_3^2)(\hat{f}_{11} + \hat{f}_{33}) - (y_1^2\hat{f}_{33} - 2y_1y_3\hat{f}_{13} + 
y_3^2\hat{f}_{11}) = \mu_{2}, $$ $$
(y_1^2 + y_2^2)(\hat{f}_{11} + \hat{f}_{22}) - (y_1^2\hat{f}_{22} - 2y_1y_2\hat{f}_{12} + 
y_2^2\hat{f}_{11}) = \mu_{3}.
\label{eqn:SysNu1}
\end{equation}
Rearranging the above gives the following
\begin{equation}
\left(\begin{array}{ccc} 0 & y_2^2 & y_3^2 \\ y_1^2 & 0 & y_3^2 \\ 
y_1^2 & y_2^2 & 0 \end{array} \right)
\left(\begin{array}{c} \hat{f}_{11}\\ \hat{f}_{22}\\ \hat{f}_{33} \end{array} \right) =  
\left(\begin{array}{c} \mu_{1} - 2y_2y_3\hat{f}_{23} \\ \mu_{2} - 2y_1y_3\hat{f}_{13} \\ \mu_{3} - 2y_1y_2\hat{f}_{12}
\end{array} \right).  
\label{eqn:SysNu2}
\end{equation}
Let us relabel the RHS of the above as
\begin{equation}
\left(\begin{array}{c} \mu_{1} - 2y_2y_3\hat{f}_{23} \\ \mu_{2} - 2y_1y_3\hat{f}_{13} \\ \mu_{3} - 2y_1y_2\hat{f}_{12}
\end{array} \right) = \left(\begin{array}{c} \nu_{1} \\ \nu_{2} \\ \nu_{3} \end{array} \right).
\label{eqn:SysNu3}
\end{equation}
In this way the solution of (\ref{eqn:SysNu2}) can be written as
\begin{equation}
\hat{f}_{11} = \frac{1}{2y_1^2}(\nu_2 + \nu_3 - \nu_1),$$ $$
\hat{f}_{22} = \frac{1}{2y_2^2}(\nu_1 + \nu_3 - \nu_2),$$ $$
\hat{f}_{33} = \frac{1}{2y_3^2}(\nu_2 + \nu_1 - \nu_3).
\label{eqn:SolnNu}
\end{equation}
Upon substitution of the off-diagonals and $\mu\text{'s}$ into (\ref{eqn:SolnNu}), we obtain
\begin{equation}
\hat{f}_{11} = \frac{1}{2y_1^2}(\mu_2 + \mu_3 - \mu_1 + y_2\lambda_2 + y_3\lambda_3 - 3y_1\lambda_1),$$ $$
\hat{f}_{22} = \frac{1}{2y_2^2}(\mu_1 + \mu_3 - \mu_2 + y_1\lambda_1 + y_3\lambda_3 - 3y_2\lambda_2),$$ $$
\hat{f}_{33} = \frac{1}{2y_3^2}(\mu_2 + \mu_1 - \mu_3 + y_2\lambda_2 + y_1\lambda_1 - 3y_3\lambda_3).
\label{eqn:SolnMu}
\end{equation}

\section{Insufficiency for two axes}

To reduce the data acquisition time, experimentalists would want to rotate the specimen on its axis as few times as possible. We show that in the general case two orthogonal axes are insufficient by considering 
components in the null space of the TRT. Thus $J_{\eta}f = 0.$ If we had two orthogonal axes, say $ \eta = e_1, e_2,$ then 
$\langle \eta, J_{\eta}f\eta \rangle = 0.$ This implies that $f_{11} = f_{22} = 0.$ From the definition of $\lambda_{\eta}(y)$ and 
$\mu_{\eta}(y),$ $\lambda_1 = \lambda_2 = 0$ and $\mu_1 = \mu_2 = 0$. The system of equations for the off-diagonals 
(\ref{eqn:SysLambda}) gives us
\begin{eqnarray}
y_2 \hat{f}_{12} + y_3 \hat{f}_{13} = \lambda_1 = 0, \\
y_1 \hat{f}_{12} + y_3 \hat{f}_{23} = \lambda_2 = 0.
\label{eqn:Lambda_2Axes}
\end{eqnarray}
Moreover the system of equations corresponding to the other \textit{non-axial} components, (\ref{eqn:SysNu3}), gives
\begin{eqnarray}
y_2^2 \hat{f}_{22} + y_3^2 \hat{f}_{33} = \mu_1 - 2y_2y_3 \hat{f}_{23}, \\
y_1^2 \hat{f}_{11} + y_3^2 \hat{f}_{33} = \mu_2 - 2y_1y_3 \hat{f}_{13}.
\label{eqn:Nu_2Axes_a}
\end{eqnarray}
Due to the argument at the start we can rearrange (\ref{eqn:Nu_2Axes_a}) as
\begin{eqnarray}
y_3^2 \hat{f}_{33} = - 2y_2y_3 \hat{f}_{23}, \\ y_3^2 \hat{f}_{33} = - 2y_1y_3 \hat{f}_{13}.
\label{eqn:Nu_2Axes_b}
\end{eqnarray}
From the above, say $\hat{f}_{33}$ is arbitrary and consequently $\hat{f}_{13} \ \text{and} \ \hat{f}_{23}$ are determined as
\begin{equation}
\hat{f}_{23} = - \frac{y_3}{2y_2} \hat{f}_{33} \ \text{and} \ \hat{f}_{13} = - \frac{y_3}{2y_1} \hat{f}_{33}.
\label{eqn:Nu_2Axes_c}
\end{equation}
Using the values obtained in (\ref{eqn:Nu_2Axes_c}) and substituting into (\ref{eqn:Lambda_2Axes}) we can write $\hat{f}_{12}$ as
\begin{eqnarray}
\hat{f}_{12} = \frac{y_3^2}{2y_1y_2} \hat{f}_{33}.
\label{eqn:Nu_2Axes_d}
\end{eqnarray}
Thus all the off-diagonal components in the tensor field are determined through $\hat{f}_{33}$ which is arbitrary. Hence two
axes are insufficient.

\section{Potential cases}

In the potential case $f_{ij}= \partial u_i /\partial x_j +\partial u_j /\partial x_i$ where $u\in \mathcal{S}(\mathbb{R}^3; \mathbb{C}^3)$. This is important for applications in that a linear strain tensor $f$ has this form where $u$ is twice the displacement field.

Without loss of generality suppose that data is known only for rotations about $\eta = e_1, e_2$. We have immediately $f_{11},f_{22}$ and hence by direct integration twice $u_1$ and $u_2$. We now also have $f_{12}$ from the partial derivatives of $u_1$ and $u_2$. It remains only to find $u_3$. 
Multiplying $\hat{f}_{12}$ by $y_1 y_2$ and $\hat{f}_{13}$ by $y_1 y_3$ and adding both of them in (\ref{eqn:SolnLambda}) gives
\begin{equation}
y_1 y_2 \hat{f}_{12} + y_1 y_3 \hat{f}_{13} = y_1 \lambda_1.
\end{equation}
This gives us $f_{13}$ in terms of known data and as $\partial u_1/\partial x_3$ is known we have $\partial u_3/\partial x_1$ and hence $u_3$. 
We summarise in the theorem
\begin{theorem}
A potential $f\in \mathcal{S}(\mathbb{R}^3; S^2(\mathbb{C}^3))$ is determined uniquely from $Jf(\xi,x)$ restricted to $\xi \in \eta_1^\perp$ and $\xi \in \eta_2^\perp$ where $\eta_1$ and $\eta_2$ are orthogonal.
\end{theorem}
This result is of considerable practical importance as it means that stain tensors, in a scheme such as that envisaged in \cite{LW}, can be recovered from rotations about only two axes.

We now show that in general a potential $f$ cannot be recovered uniquely from a one-axis rotation by constructing a general element of the null space.  
Suppose we rotate only about $e_1$ we have immediately $f_{11}=0$ and as $\hat{f}_{ij} = y_i \hat{u}_j + y_j \hat{u}_i$ we see $u_1=0$. Now as $\lambda_1$ and $u_1$ are zero
\begin{equation}
y_2  \hat{u}_2 + y_3  \hat{u}_3 = \frac{\left(\lambda_1 - (y_2^2 +y_3^2)\hat{u}_1\right)}{y_1} = 0. 
\end{equation}
and as $\mu_1=0$
\begin{equation}
 y_2\hat{u}_2 + y_3\hat{u}_3 = \frac{\mu_1}{2(y_2^2+ y_3^2)}=0,
\end{equation}
giving no new information. So $u$ must satisfy $u_1=0$ and $\partial u_2/\partial x_2 + \partial u_3/\partial x_3 =0$. 
For example if $u_2$ is arbitrarily specified, then
\begin{equation}
u_3 = -\int\limits_{-\infty}^{\infty} \frac{\partial u_2}{\partial x_2}\, \mathrm{d}x_3.
\end{equation}

\section{Numerical results}

\subsection{Forward model}

In order to generate data sets, we need to implement a discretized version of the operator $J$ described in (\ref{eqn:TRT}) as 
a matrix which will calculate integrals of projections and act upon generated strain fields represented by an array. Instead of 
calculating the whole matrix, we generate it one row at a time (on the {\em fly}) which corresponds to one individual 
source-detector pair for one of the three components in $P_{\xi}f$ described by (\ref{eqn:J_project}).  

\subsubsection{Discrete representation of the tensor field}

The discretized tensor  field is stored as a $ N \times N \times N \times 6 $ vector, containing the 6 distinct 
values of the symmetric second rank tensor field for each voxel in a $ N \times N \times N $ voxel grid. We increment first by the 
tensor component number, then the position $x_1$, $x_2$ and finally $x_3$. Furthermore, the data is represented by a
$ 3 \times n_{\theta} \times h \times w \times 3 $ multi-dimensional array (5 dimensional), where we use $3$ rotation axes 
$(\eta = e_1,e_2 \ \text{and} \ e_3)$, $n_{\theta}$ angles steps for tomographic acquisition around each axis and $ h \times w $
represents how many rays in the horizontal and vertical direction. The factor of $3$ is the number of independent values of $Jf$ 
in (\ref{eqn:J_project}) which we integrate along each ray. 

\subsubsection{Methodology}

We simulate an experimental setup with parallel rays passing through a specimen. Sources and detectors 
consist of arrays in an equally spaced grid, either side of the object being scanned. The ratio of the number of rays in the 
horizontal to vertical direction is $4:3$. The source-detector pair is kept fixed and the object is rotated. 
We follow the procedure of   \cite[Chapter 5.1.4]{DS} to calculate the approximate integral along a line through a voxel grid. This will
give us the contribution of each voxel to the total integral for a given ray. For a given tensor we need to calculate the projection 
on to the plane perpendicular to the ray. To simplify this, we rotate the coordinate system. We extract
the two appropriate components (relating to the axial and non-axial components). Since we have the 
contribution of each voxel to the integral, the length of intersection of the ray with the voxel,  from ray tracing we can form the sum of these intersection lengths with the voxel values to form the approximate integral.

For our numerical experiments, phantoms were generated inside a cubic grid measuring $405 \times 405 \times 405$ voxels and 
measurements were simulated for a source/detector grid with $ 405 \times 540 $ pixels. The number of rays in the vertical direction was taken
to be the image height (i.e. 405 pixels). The pixel grid was then down-sampled by a factor of 3 to $135 \times 180$ pixels by
the process of binning. Finally $1\%$ Gaussian pseudo-random noise was added before reconstructing on a courser voxel grid 
measuring $135 \times 135 \times 135$. The number of angles (projections) was $240$ per rotation axis. 

\subsection{Generating Phantoms}

For input into the forward projector we generate two different phantoms or test fields. The first one only has smooth features 
and is expected to be less sensitive to algorithmic instabilities. The second phantom contains sharp 
edges and is designed to highlight the limitations of the explicit reconstruction algorithm for discontinuous strain fields.  

\subsubsection{Phantom 1: smooth}

Phantom 1 is constructed from smooth Gaussians which makes it relatively easy to reconstruct.  We define a cubic domain $[-1,1]^3$ on which the components of $f$ are supported, defined by 
$3$-dimensional Gaussians $b_{\alpha}(x)$ for each of the components $f_{ij}$ according to Table \ref{tab:smooth_phantom}, where
$$ b_\alpha(x) = \alpha \exp(-50| x - a |^2).$$ 

\begin{table}[!htb]
    \begin{minipage}{.5\linewidth}
      \centering
      \caption{Phantom 1 - Smooth}
      \begin{tabular}{|c|r|r|r|r|} 
      \hline
      $f_{ij}$ & $\alpha$ & $a_1$ & $a_2$ & $a_3$ \\  
      \hline\hline
      \multirow{3}{*}{$f_{11}$} & -1 & -0.5 & -0.5 & -0.5 \\ 
      & 1 & -0.5 & 0.5 & -0.5 \\
      & -1 & -0.5 & 0.5 & 0.5 \\ \hline
      \multirow{2}{*}{$f_{12}$} & 1 & 0.5 & -0.5 & 0.5 \\
      & -1 & 0.5 & 0.5 & -0.5 \\ \hline
      \multirow{3}{*}{$f_{13}$} & 1 & -0.5 & -0.5 & -0.5 \\ 
      & -1 & -0.5 & -0.5 & 0.5 \\
      & 1 & -0.5 & 0.5 & 0.5 \\ \hline
      \multirow{3}{*}{$f_{22}$} & -1 & 0.5 & -0.5 & -0.5 \\ 
      & 1 & 0.5 & 0.5 & 0.5 \\
      & -1 & 0.5 & 0.5 & 0.5 \\ \hline
      \multirow{2}{*}{$f_{23}$} & 1 & -0.5 & -0.5 & 0.5 \\
      & -1 & -0.5 & 0.5 & -0.5 \\ \hline
      \multirow{3}{*}{$f_{33}$} & 1 & 0.5 & -0.5 & -0.5 \\ 
      & -1 & 0.5 & -0.5 & 0.5 \\
      & 1 & 0.5 & 0.5 & 0.5 \\ \hline
      \end{tabular}
    \label{tab:smooth_phantom}  
    \end{minipage}%
    \begin{minipage}{.5\linewidth}
      \centering
        \caption{Phantom 2 - Sharp}
        \begin{tabular}{|c|c|c|c|c|} 
      \hline
      $i$ & $j$ & $I_1$ & $I_2$ & $I_3$ \\  
      \hline\hline
      1 & 1 & [-0.4,0.4] & [-0.6,0.2] & [-0,8,0.8] \\ 
      1 & 2 & [-0.4,0.4] & [-0.2,0.6] & [-0.8,0.8] \\
      1 & 3 & [-0.8,0.8] & [-0.4,0.4] & [-0.6,0.2] \\
      2 & 2 & [-0.8,0.8] & [-0.4,0.4] & [-0.2,0.6] \\
      2 & 3 & [-0.6,0.2] & [-0.8,0.8] & [-0.4,0.4] \\ 
      3 & 3 & [-0.2,0.6] & [-0.8,0.8] & [-0.4,0.4] \\
      \hline
      \end{tabular}
    \label{tab:sharp_phantom}  
    \end{minipage}
\end{table}

\subsubsection{Phantom 2: sharp}

Phantom 2 is constructed to contain sharp edges, to highlight non-linear strain fields which is not quite compatible with the
reconstruction algorithm. As in the smooth case, we define $f$ on the cube $[-1,1]^3$, but set $f_{ij}$ to be the characteristic function of $I_1\times I_2 \times I_3$ , according to Table \ref{tab:sharp_phantom}. 

\subsection{Reconstruction Procedure}

The recovery of {\em axial} components are relatively straight forward as this is just plane-by-plane Radon inversion. 
Hence, we apply a ramp-filter (Ram-Lak) to the data and back-project to achieve the diagonal entries for each rotation axis.
Since, we can think of back-projection as the dual operator of ray integration, we reuse the ray tracing code to implement  
back-projection as the transpose of ray integration. 

From (\ref{eqn:Lemma4.1b}), we see that the simulated data values $J^1$ that are collected for each plane need to be 
differentiated in the $p$-direction, before any back-projection takes place. To implement this, we perform a regularised 
derivative, which we carry out in the Fourier domain using a Hamming window regularisation. After performing a one-dimensional 
Discrete Fourier Transform  using the Fast Fourier Transform FFT algorithm, we multiply by $ -ikw(n)$, where $ w(n)$ is the Hamming window.  For a discrete signal of length N, labelled by $n= 0,....,N-1$, the {\em Hamming window} $w(n)$ is given by
\begin{equation}
 w(n) = 0.54 - 0.46\cos\left(\frac{2\pi n}{N-1}\right).
\label{eqn:Hamming_window} 
\end{equation}
The result is returned to the spatial domain using the inverse FFT algorithm.

Following the filter we back-project the differentiated plane by plane data onto the voxel grid and the tangential vector field 
components (i.e. $\lambda$) are calculated using a three dimensional FFT algorithm and application of a {\em ramp-filter} in 
frequency space. Then equations (\ref{eqn:SolnLambda}) are used to recover the off-diagonal terms in frequency space. The only 
exception is the voxel $(y_1,y_2,y_3) = (0,0,0)$, where $\hat{f}_{12},\hat{f}_{13}$ and $\hat{f}_{23}$ are undefined. Here, 
the value is set using linear interpolation from nearby voxels. To complete our reconstruction, we employ the three 
dimensional inverse FFT to recover $f_{ij}$.

\subsection{Results and summary}

Below we illustrate the results of our implemented reconstructions which clearly show the performance of the algorithm on two 
different tensor fields. As expected, the smooth  (Gaussian) field is reconstructed well. However, when we 
introduce sharp edges in the components of a field such as a crack, we see that as expected the 
reconstruction is inaccurate and artefacts are visible. The artefacts increase as more noise is added.   
%\graphicspath{{./numeric_results/}}

\begin{figure}
    \centering
    \begin{subfigure}[b]{0.45\textwidth}
        \centering
        \includegraphics[width=\textwidth]{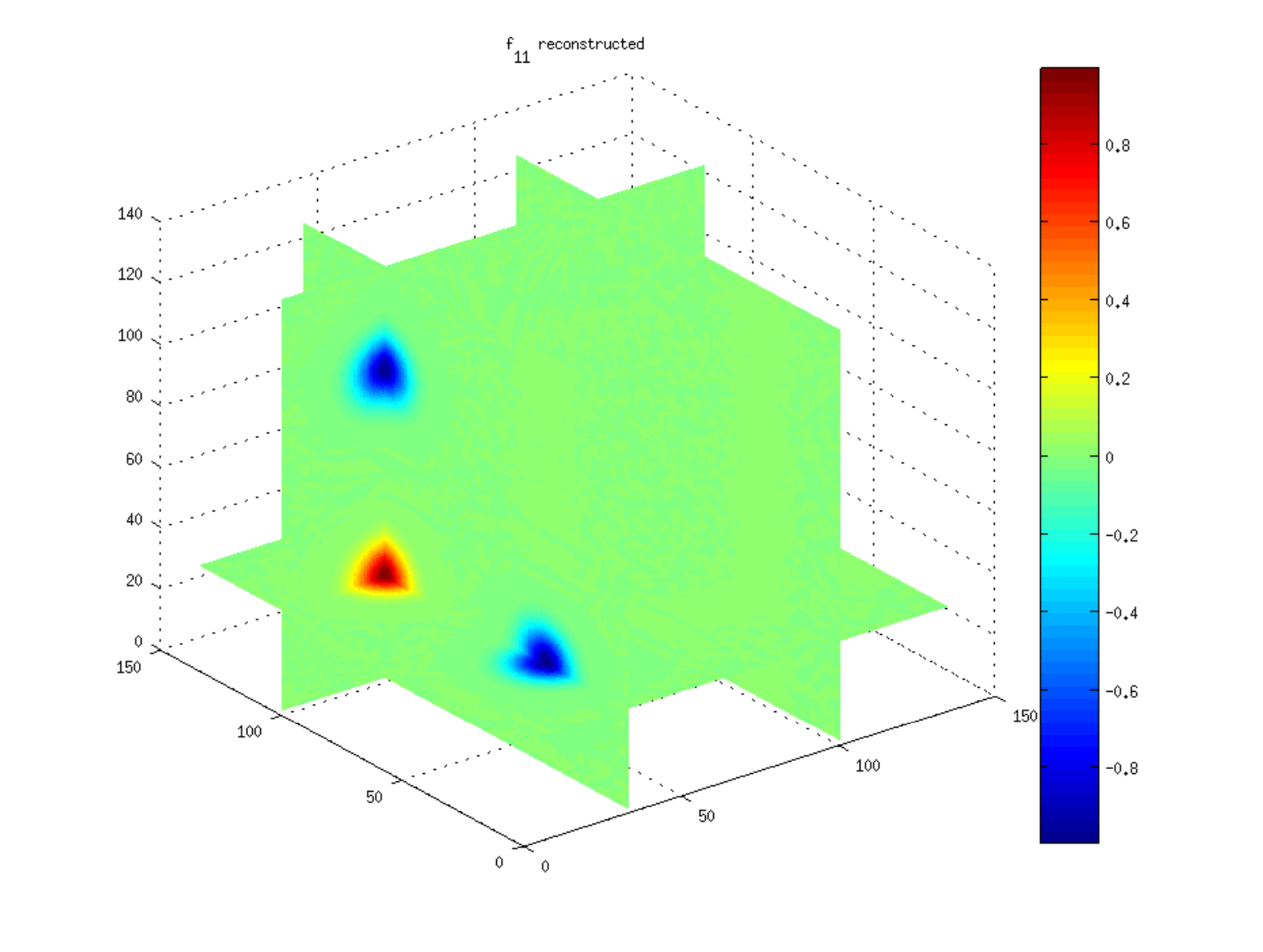}
        \caption{Reconstruction of $f_{11}$}
	\label{fig:rec_smooth2_f11}
    \end{subfigure}
    \hspace{5mm}
    \begin{subfigure}[b]{0.45\textwidth}
        \centering
        \includegraphics[width=\textwidth]{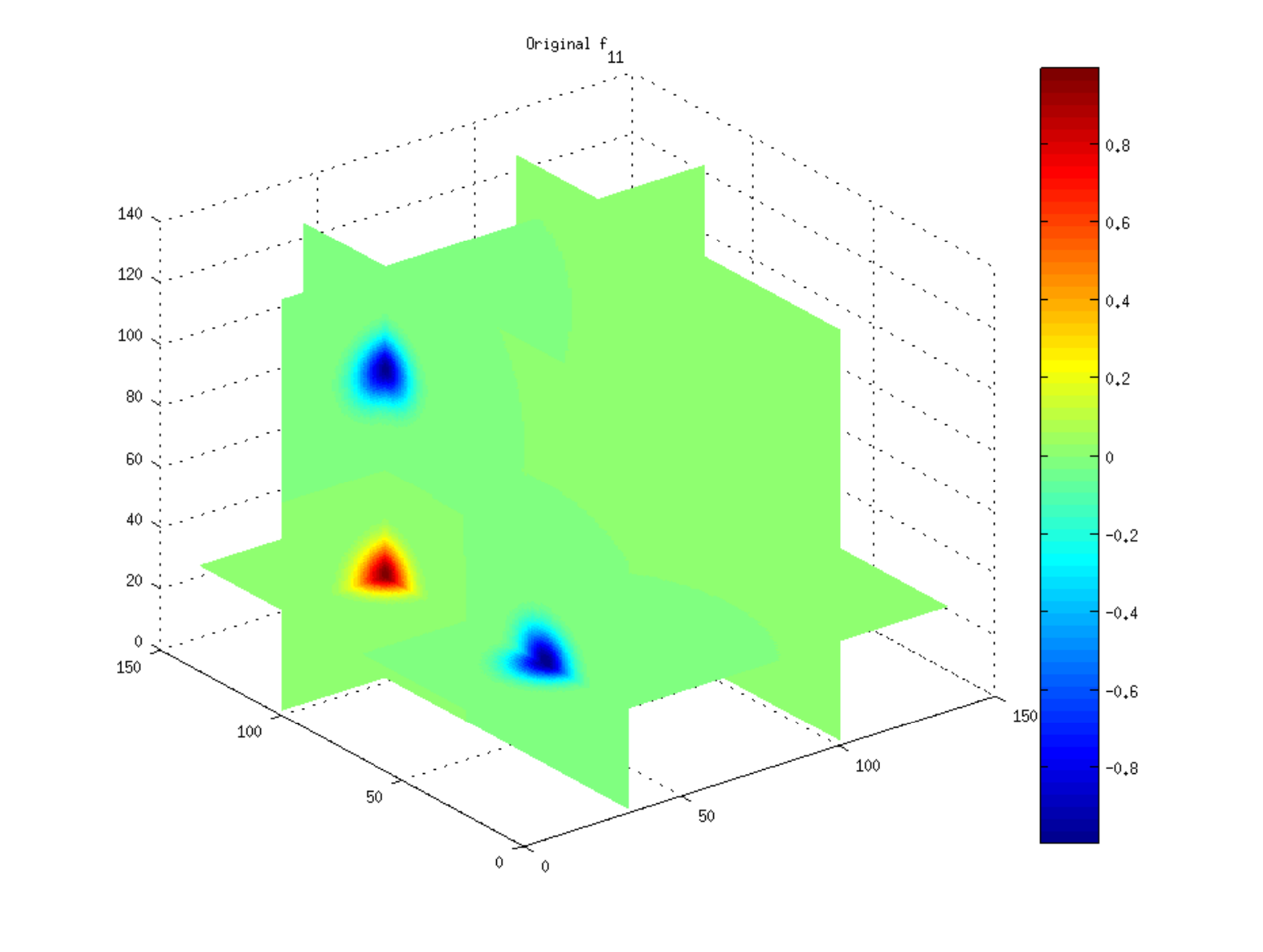}
        \caption{Original $f_{11}$}
	\label{fig:ori_smooth2_f11}
    \end{subfigure}
    \vspace{5mm}
    \begin{subfigure}[b]{0.45\textwidth}
        \centering
        \includegraphics[width=\textwidth]{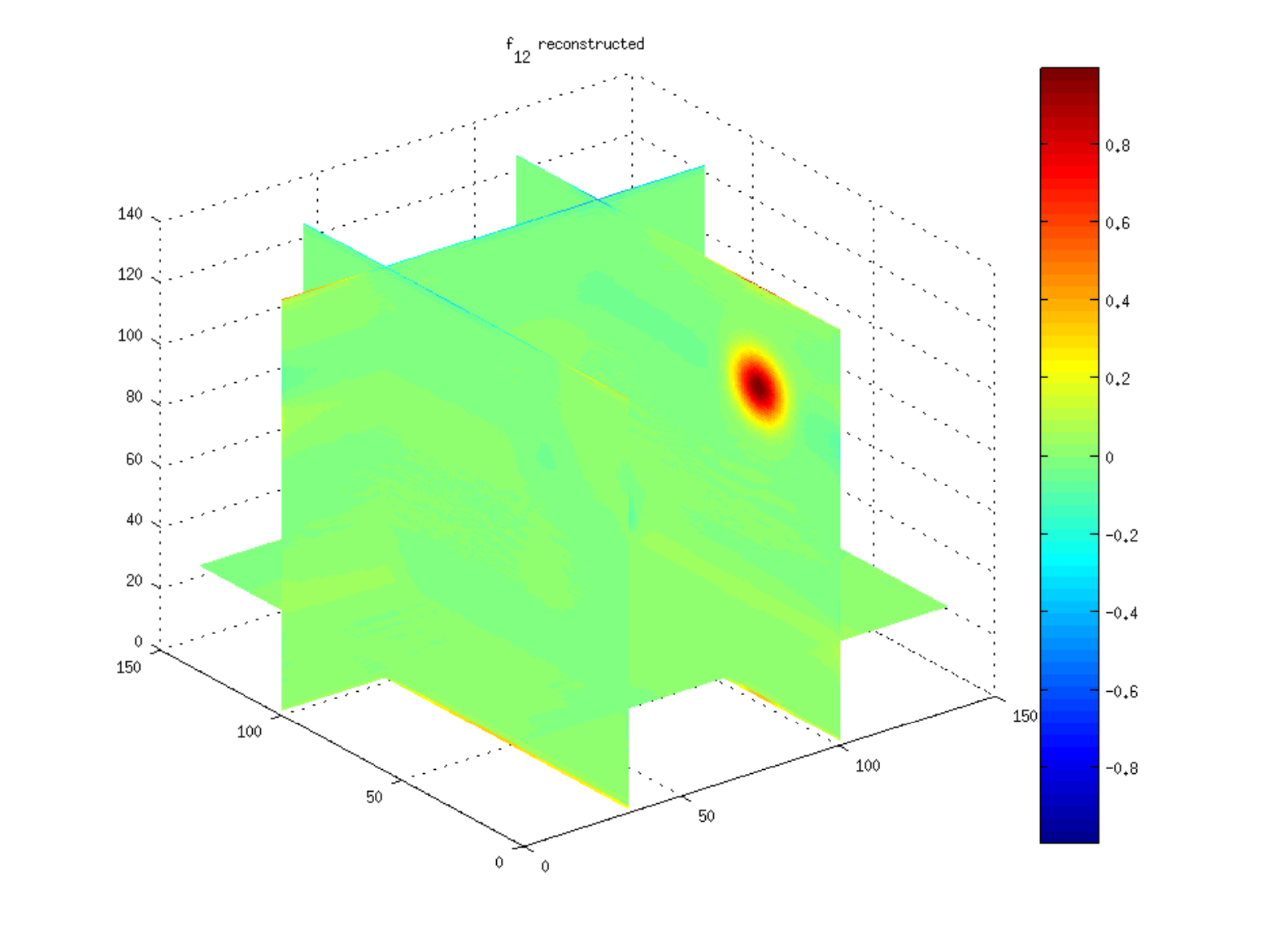}
        \caption{Reconstruction of $f_{12}$}
        \label{fig:rec_smooth2_f12}
    \end{subfigure}
    \hspace{5mm}
    \begin{subfigure}[b]{0.45\textwidth}
        \centering
        \includegraphics[width=\textwidth]{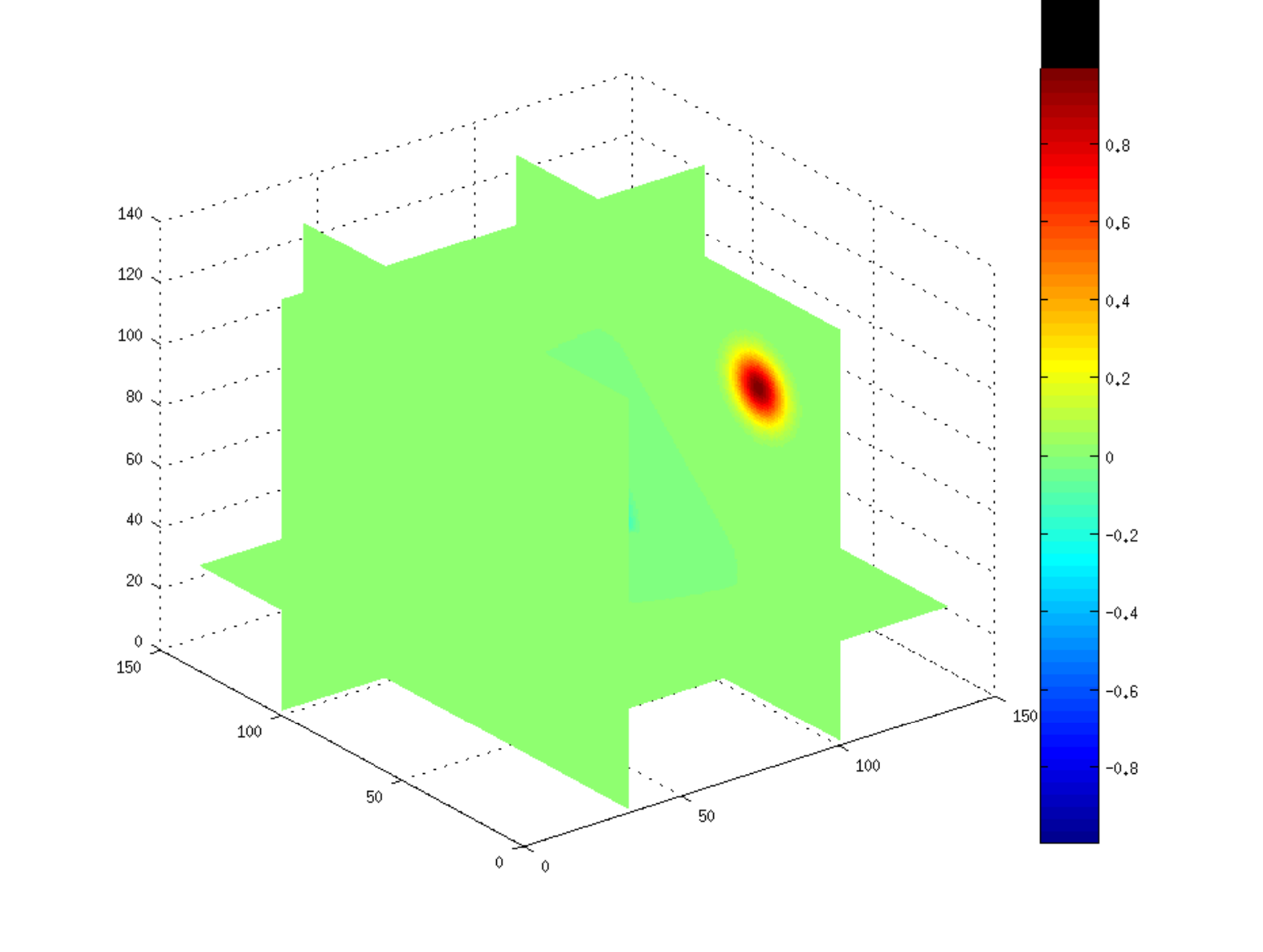}
        \caption{Original $f_{12}$}
        \label{fig:ori_smooth2_f12}
    \end{subfigure}
    \vspace{5mm}
    \begin{subfigure}[b]{0.45\textwidth}
        \centering
        \includegraphics[width=\textwidth]{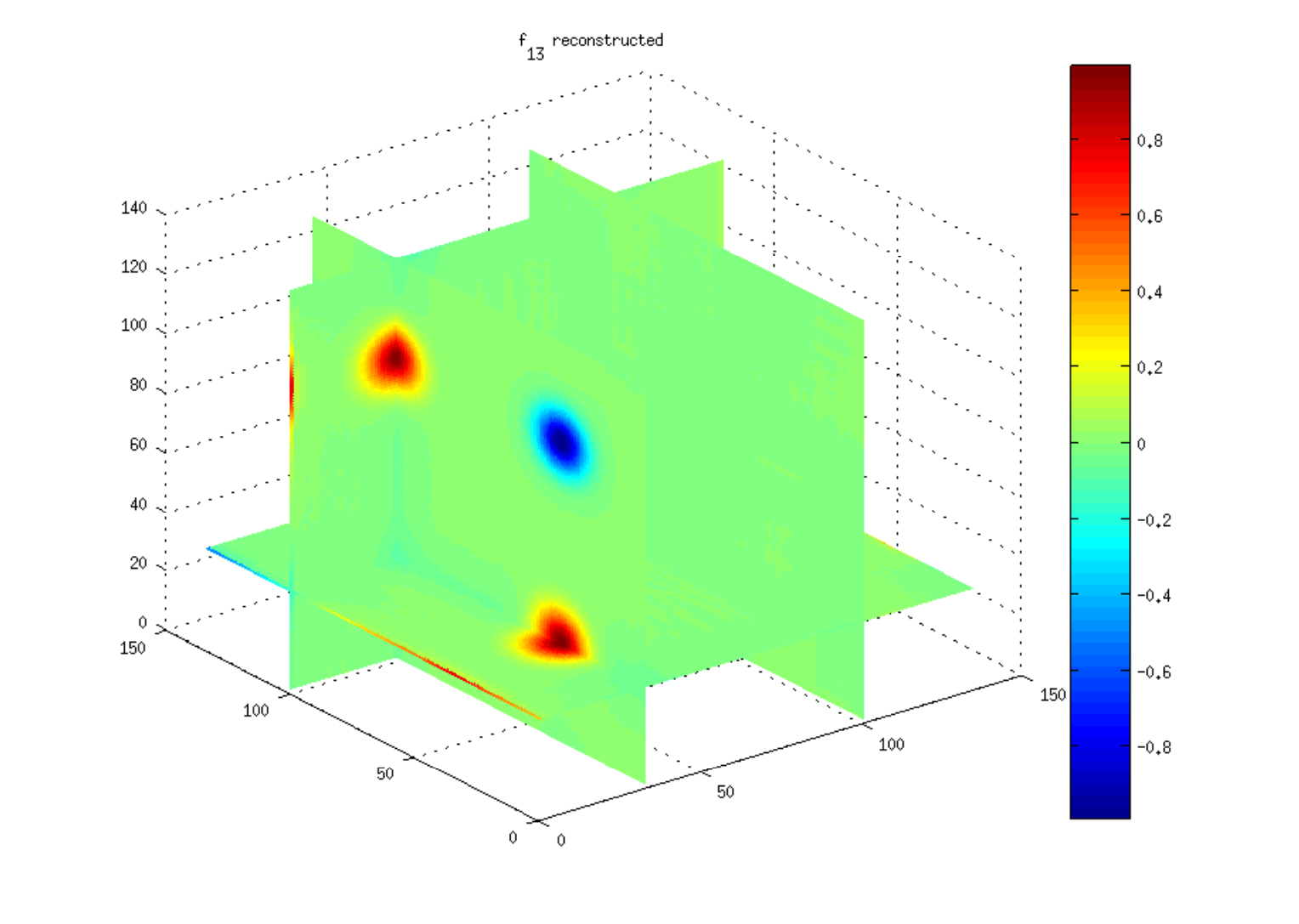}
        \caption{Reconstruction of $f_{13}$}
        \label{fig:rec_smooth2_f13}
    \end{subfigure}
    \hspace{5mm}
    \begin{subfigure}[b]{0.45\textwidth}
        \centering
        \includegraphics[width=\textwidth]{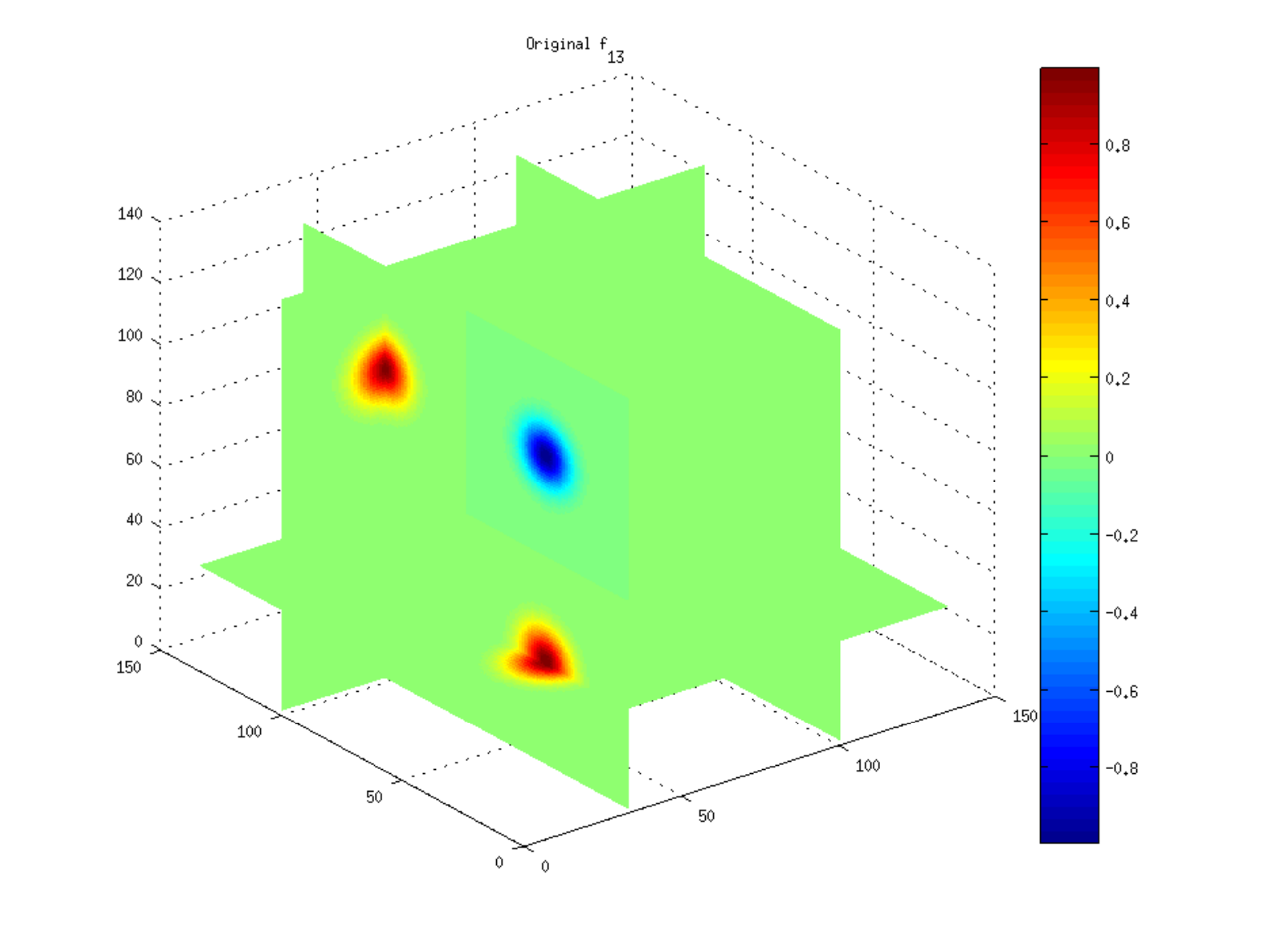}
        \caption{Original $f_{13}$}
        \label{fig:ori_smooth2_f13}
    \end{subfigure}
    \caption{Reconstruction of several gaussian balls $f_{11}, f_{12}$ and $f_{13}$}
    \label{fig:smooth2_f11_to_f13}
\end{figure}

\begin{figure}
    \centering
    \begin{subfigure}[b]{0.45\textwidth}
        \centering
        \includegraphics[width=\textwidth]{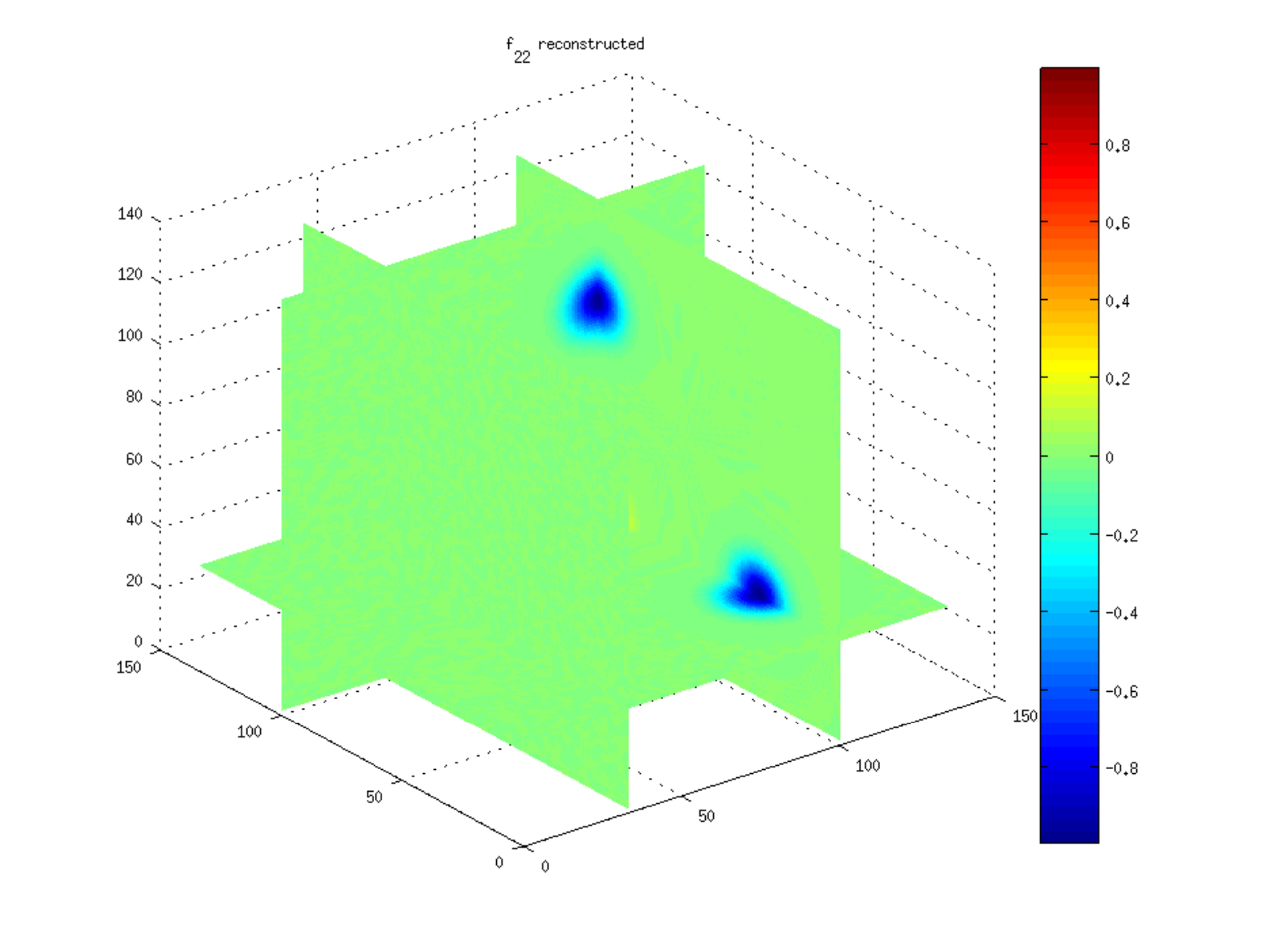}
        \caption{Reconstruction of $f_{22}$}
	\label{fig:rec_smooth2_f22}
    \end{subfigure}
    \hspace{5mm}
    \begin{subfigure}[b]{0.45\textwidth}
        \centering
        \includegraphics[width=\textwidth]{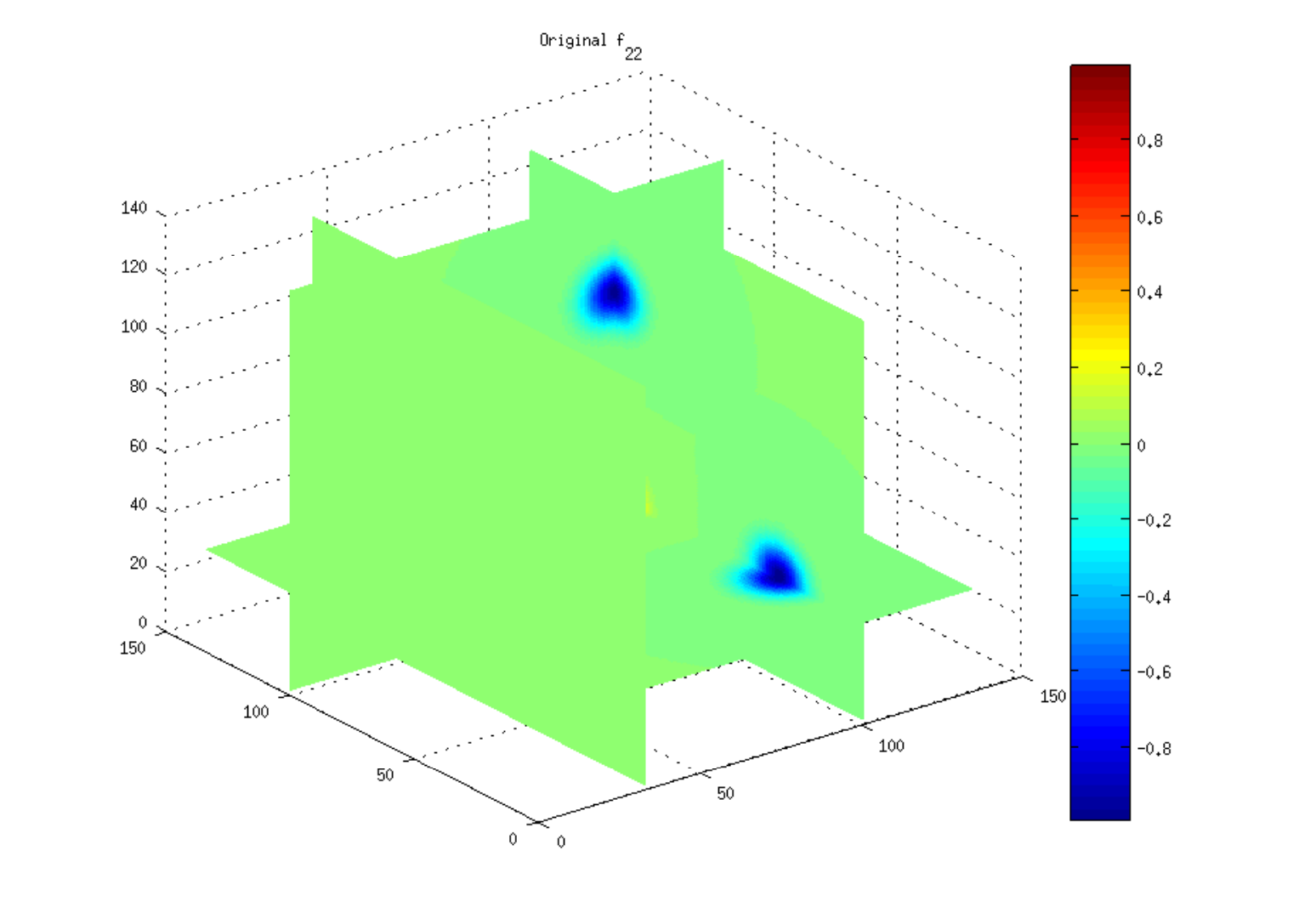}
        \caption{Original $f_{22}$}
	\label{fig:ori_smooth2_f22}
    \end{subfigure}
    \vspace{5mm}
    \begin{subfigure}[b]{0.45\textwidth}
        \centering
        \includegraphics[width=\textwidth]{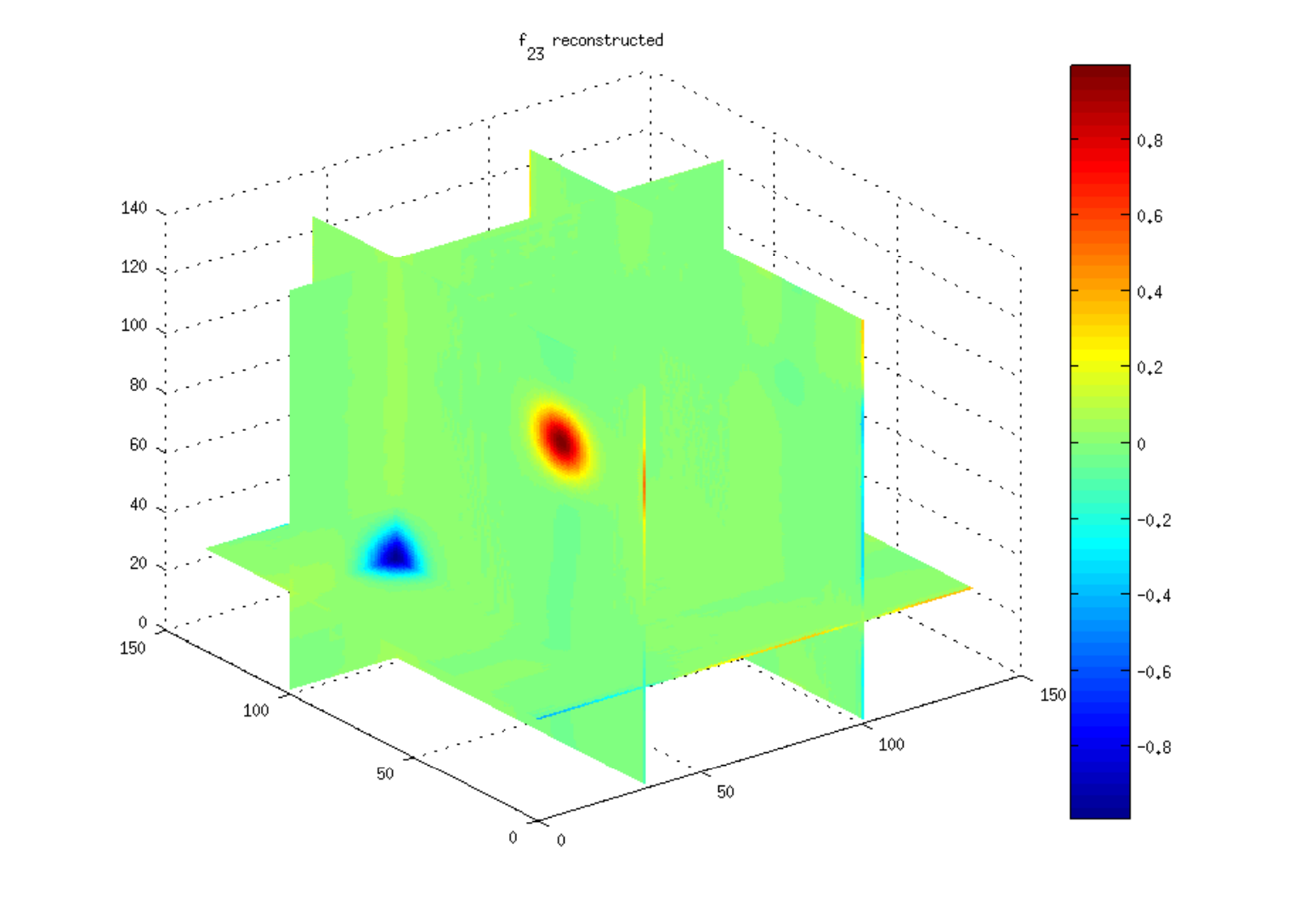}
        \caption{Reconstruction of $f_{23}$}
        \label{fig:rec_smooth2_f23}
    \end{subfigure}
    \hspace{5mm}
    \begin{subfigure}[b]{0.45\textwidth}
        \centering
        \includegraphics[width=\textwidth]{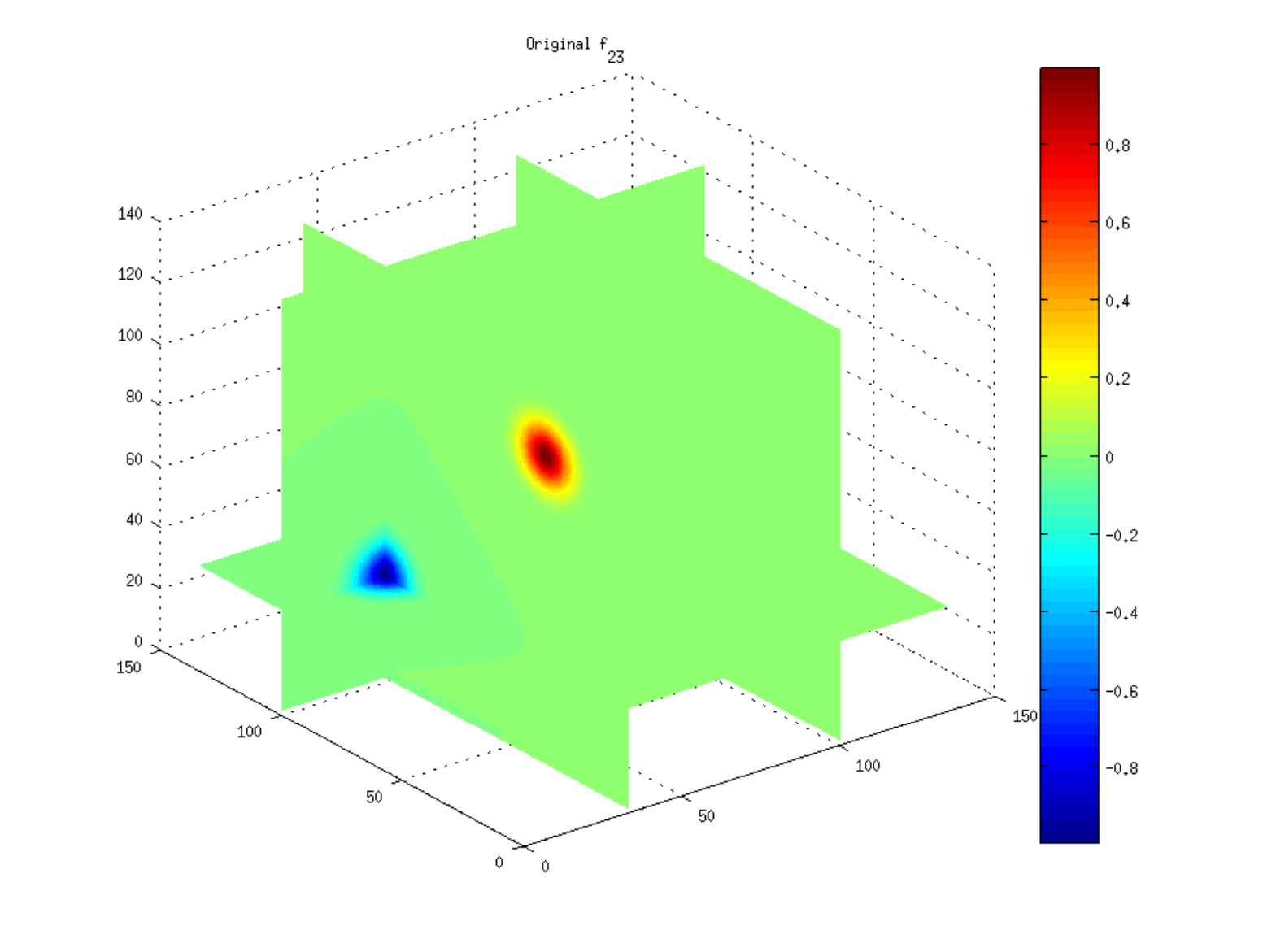}
        \caption{Original $f_{23}$}
        \label{fig:ori_smooth2_f23}
    \end{subfigure}
    \vspace{5mm}
    \begin{subfigure}[b]{0.45\textwidth}
        \centering
        \includegraphics[width=\textwidth]{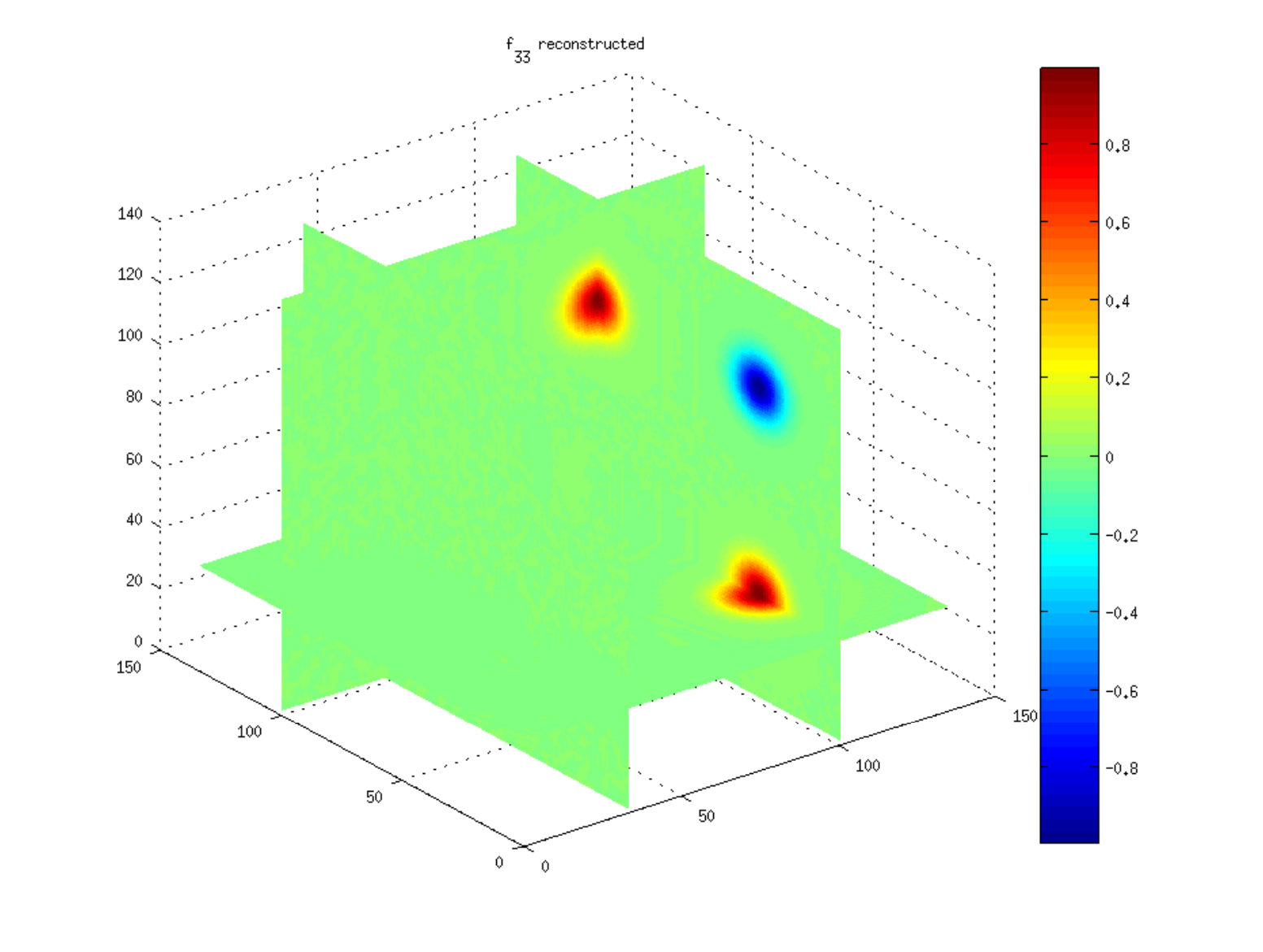}
        \caption{Reconstruction of $f_{33}$}
        \label{fig:rec_smooth2_f33}
    \end{subfigure}
    \hspace{5mm}
    \begin{subfigure}[b]{0.45\textwidth}
        \centering
        \includegraphics[width=\textwidth]{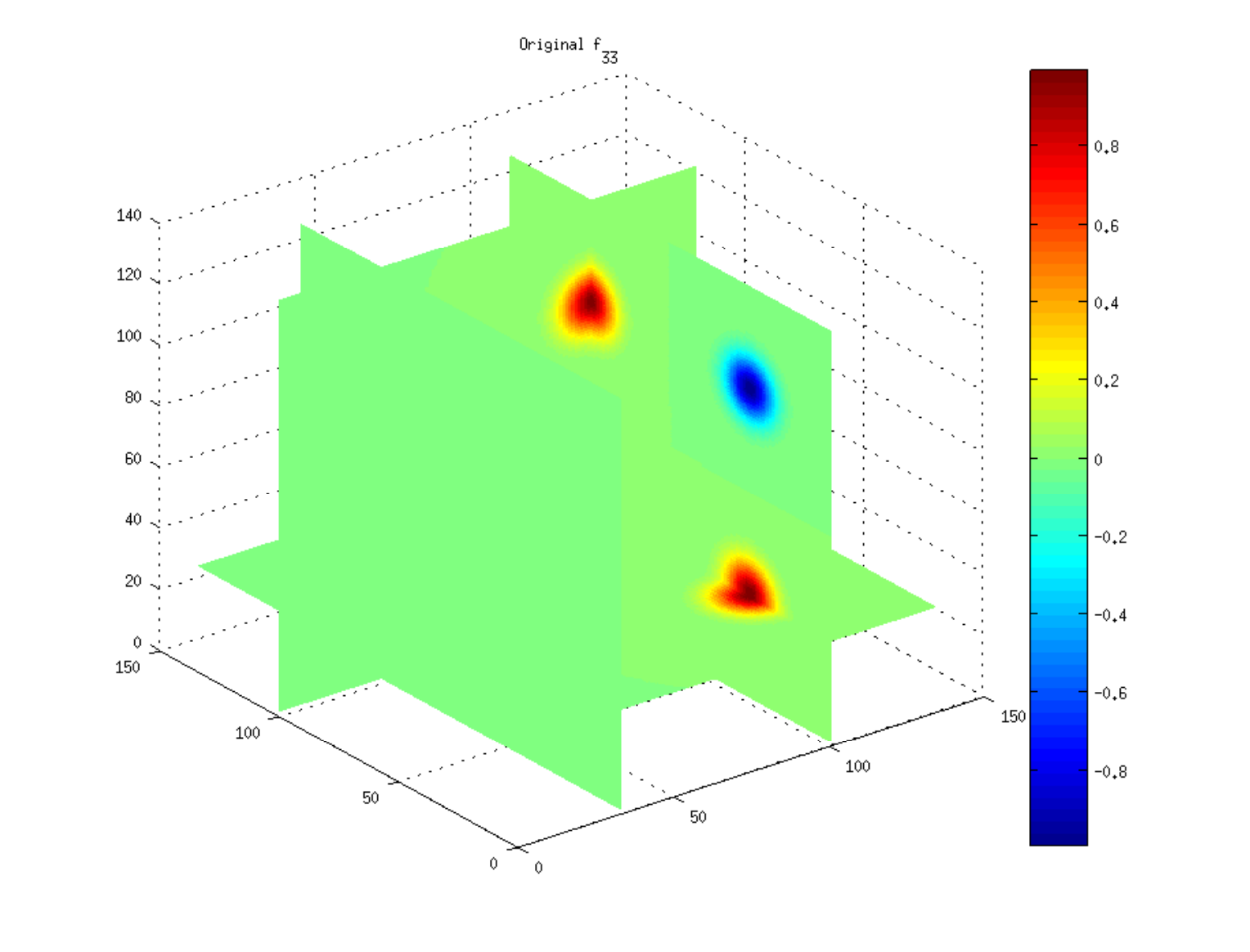}
        \caption{Original $f_{33}$}
        \label{fig:ori_smooth2_f33}
    \end{subfigure}
    \caption{Reconstruction of several gaussian balls for $f_{22}, f_{23}$ and $f_{33}$}
    \label{fig:smooth2_f22_to_f33}
\end{figure}

\begin{figure}
    \centering
    \begin{subfigure}[b]{0.45\textwidth}
        \centering
        \includegraphics[width=\textwidth]{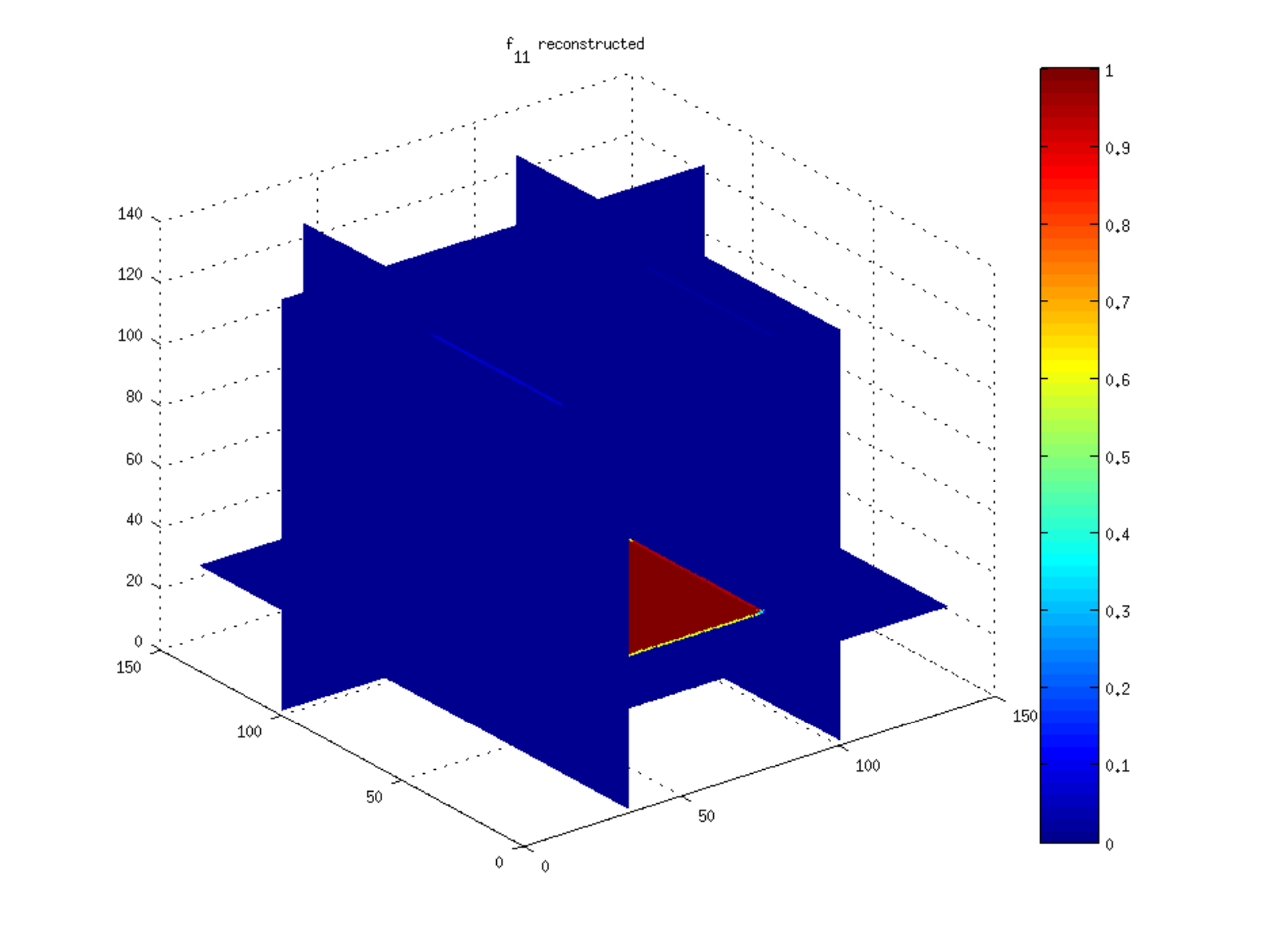}
        \caption{Reconstruction of $f_{11}$}
	\label{fig:rec_sharp_f11}
    \end{subfigure}
    \hspace{5mm}
    \begin{subfigure}[b]{0.45\textwidth}
        \centering
        \includegraphics[width=\textwidth]{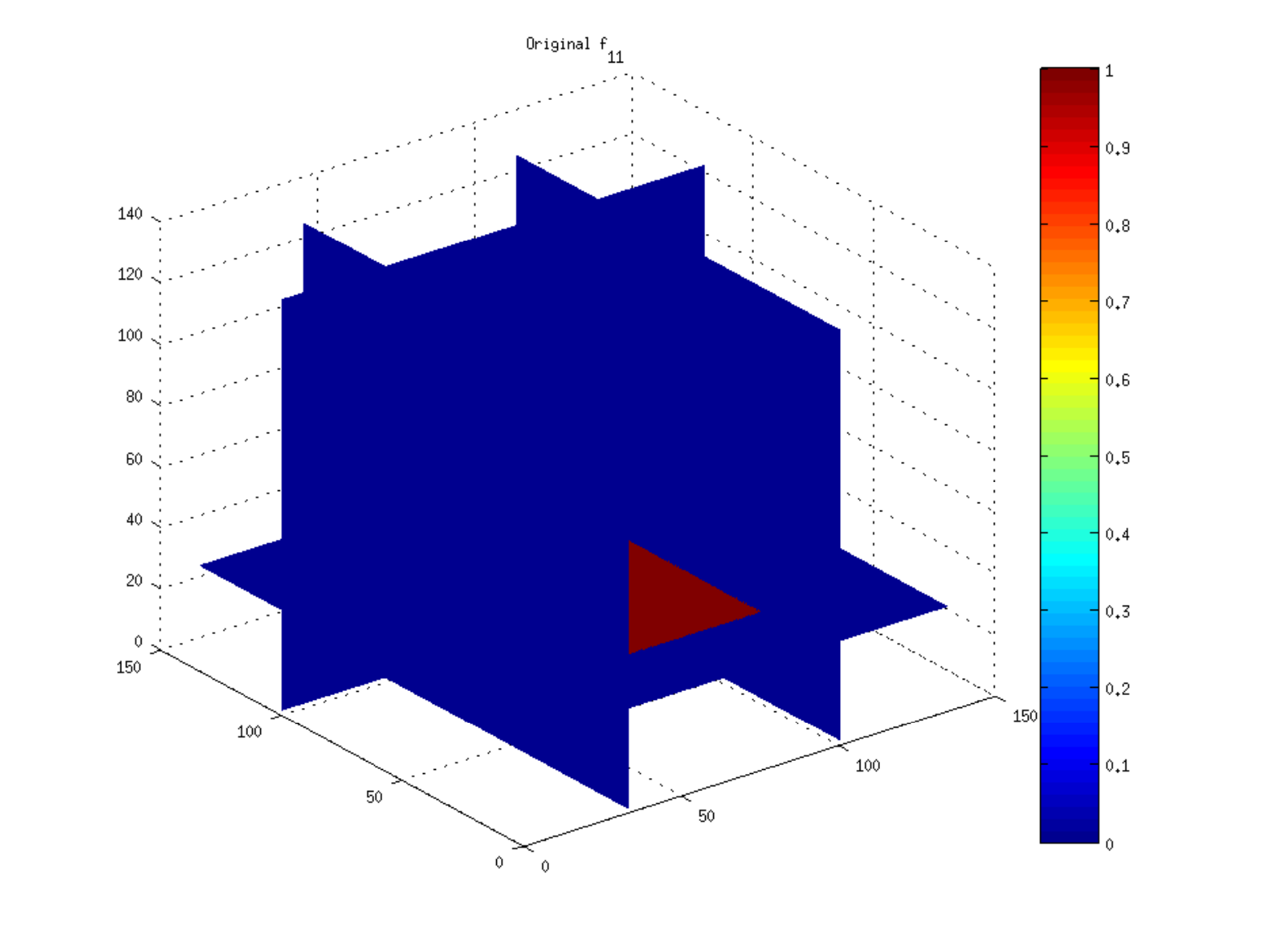}
        \caption{Original $f_{11}$}
	\label{fig:ori_sharp_f11}
    \end{subfigure}
    \vspace{5mm}
    \begin{subfigure}[b]{0.45\textwidth}
        \centering
        \includegraphics[width=\textwidth]{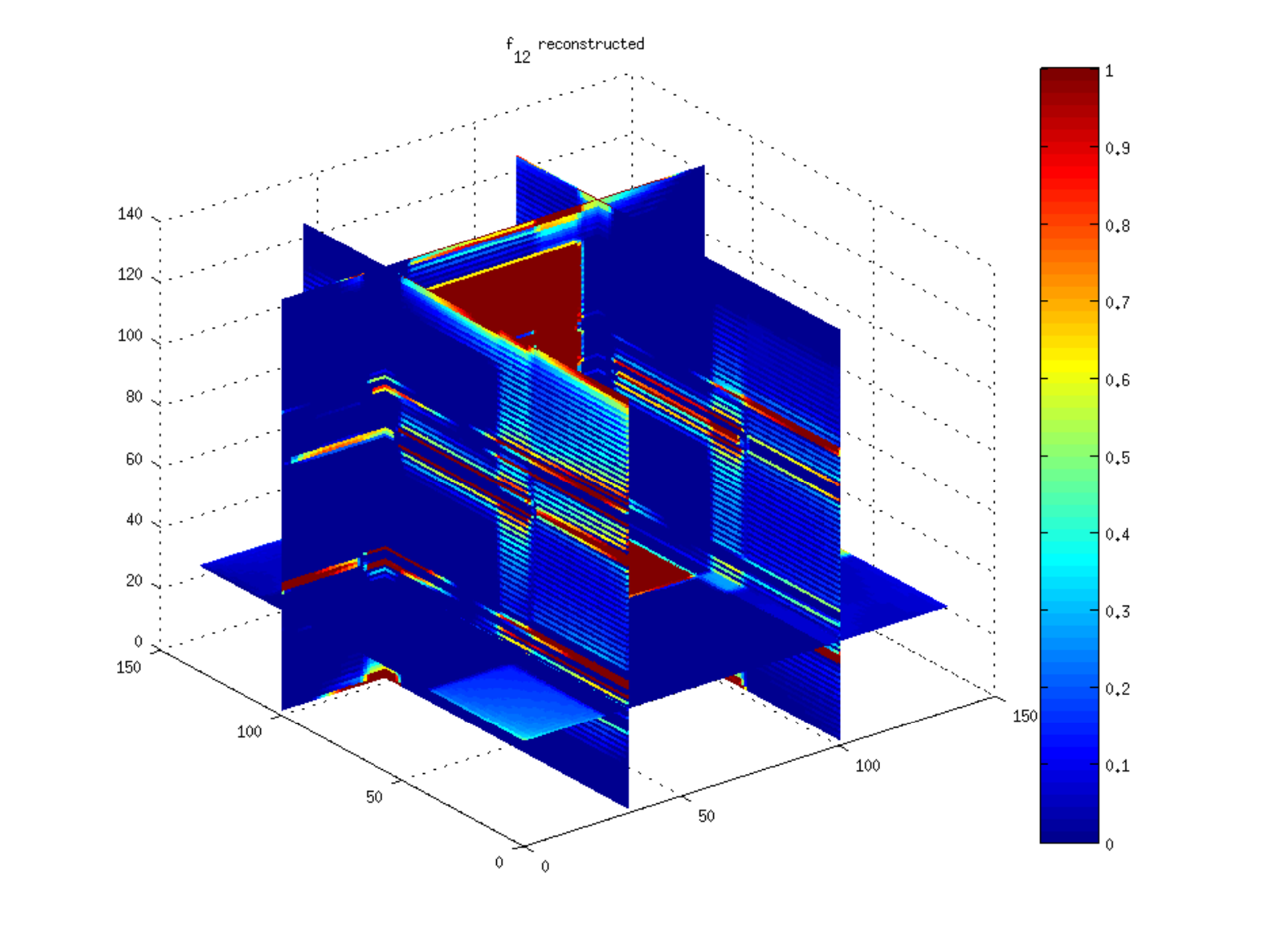}
        \caption{Reconstruction of $f_{12}$}
        \label{fig:rec_sharp_f12}
    \end{subfigure}
    \hspace{5mm}
    \begin{subfigure}[b]{0.45\textwidth}
        \centering
        \includegraphics[width=\textwidth]{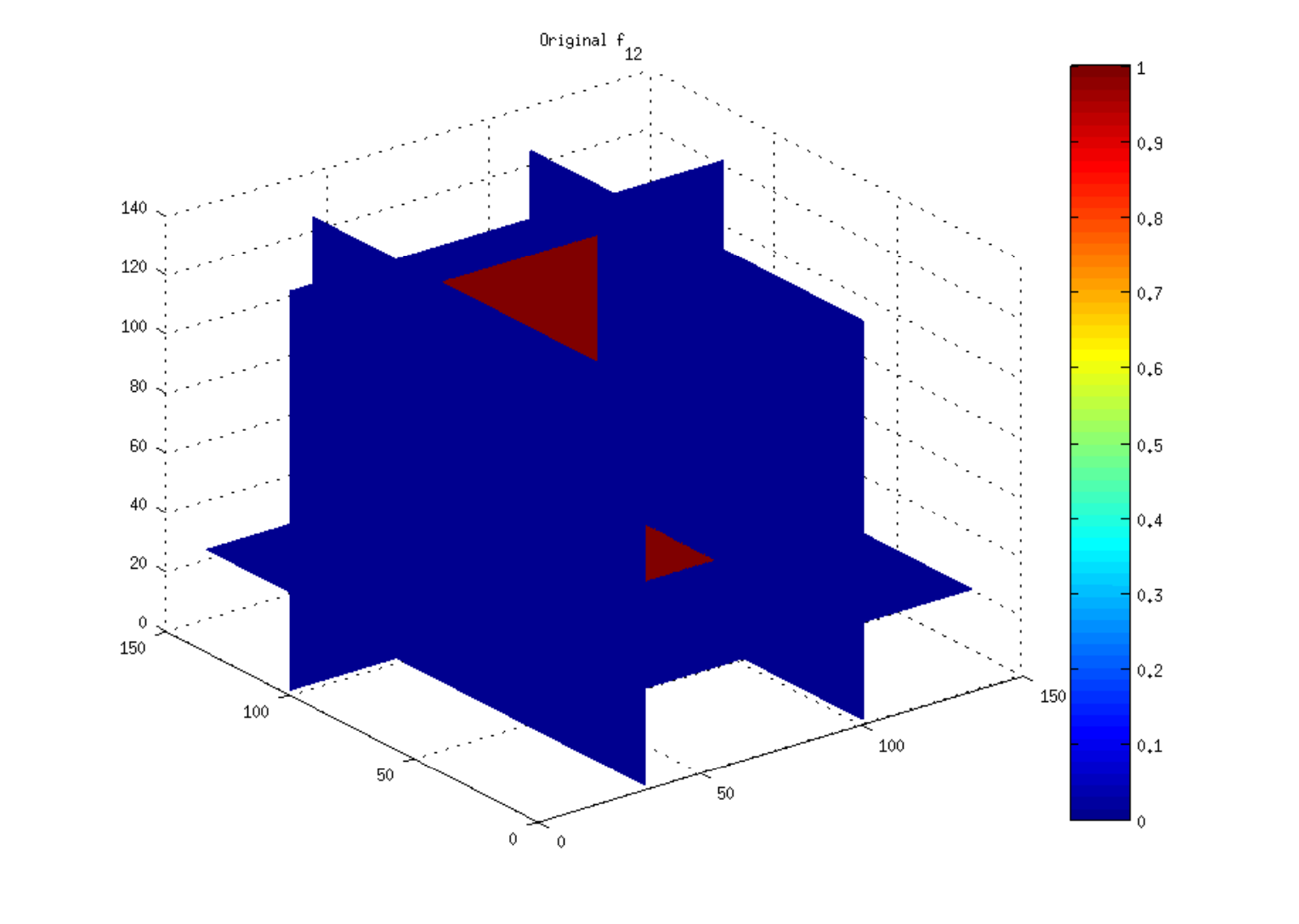}
        \caption{Original $f_{12}$}
        \label{fig:ori_sharp_f12}
    \end{subfigure}
    \vspace{5mm}
    \begin{subfigure}[b]{0.45\textwidth}
        \centering
        \includegraphics[width=\textwidth]{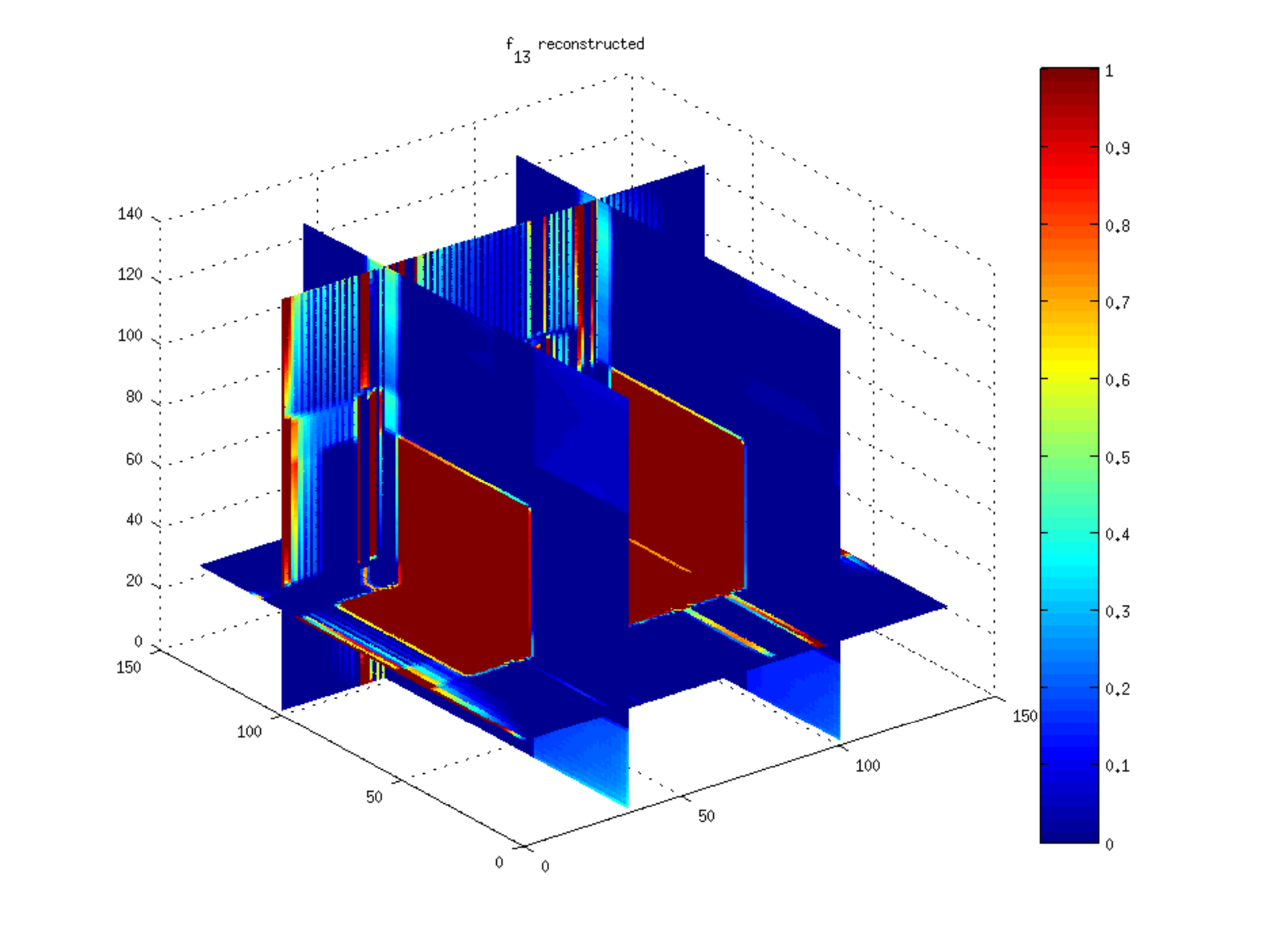}
        \caption{Reconstruction of $f_{13}$}
        \label{fig:rec_sharp_f13}
    \end{subfigure}
    \hspace{5mm}
    \begin{subfigure}[b]{0.45\textwidth}
        \centering
        \includegraphics[width=\textwidth]{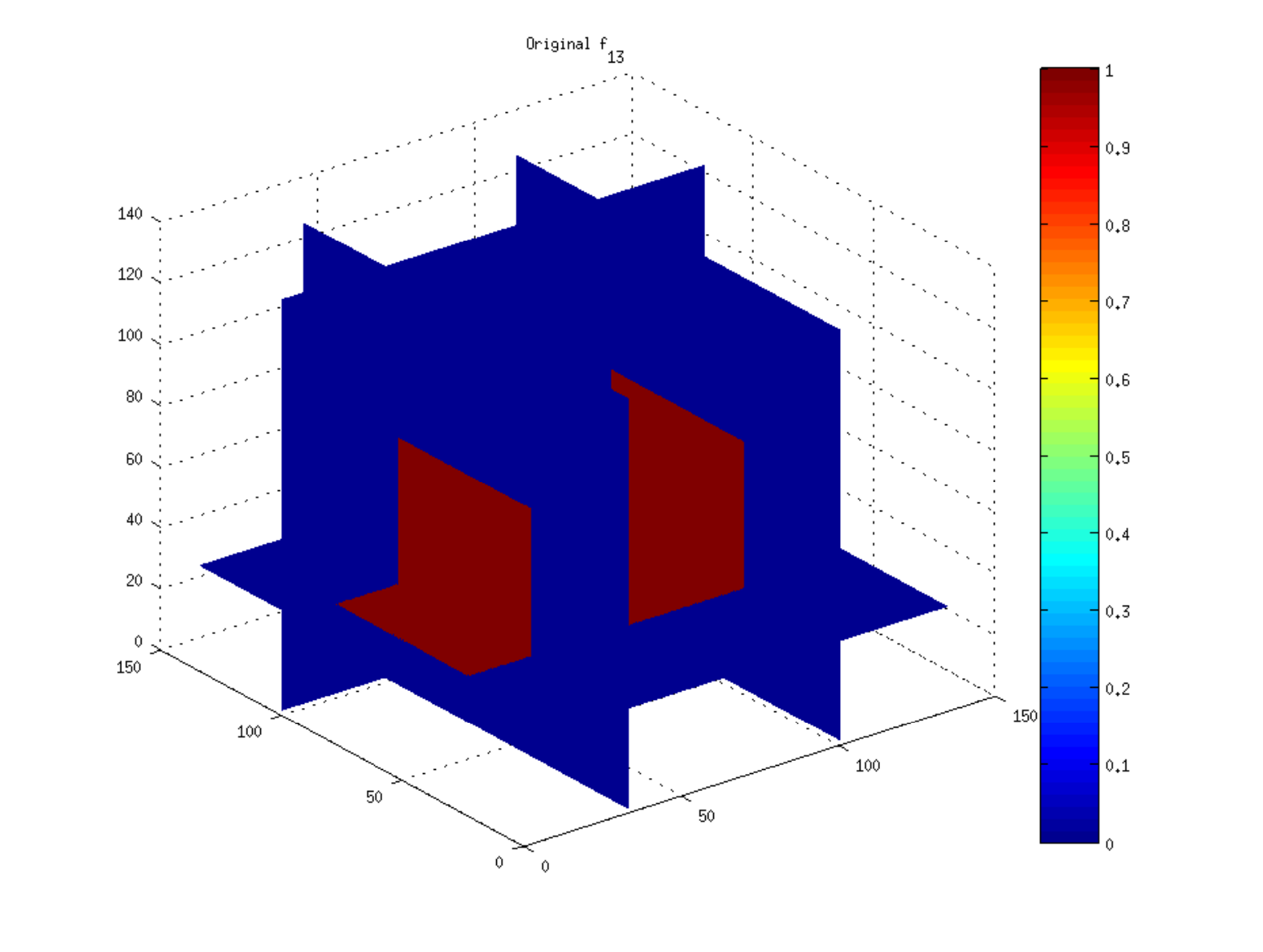}
        \caption{Original $f_{13}$}
        \label{fig:ori_sharp_f13}
    \end{subfigure}
    \caption{Reconstruction of sharp phantom, $f_{11}, f_{12}$ and $f_{13}$}
    \label{fig:sharp_f11_to_f13}
\end{figure}

\begin{figure}
    \centering
    \begin{subfigure}[b]{0.45\textwidth}
        \centering
        \includegraphics[width=\textwidth]{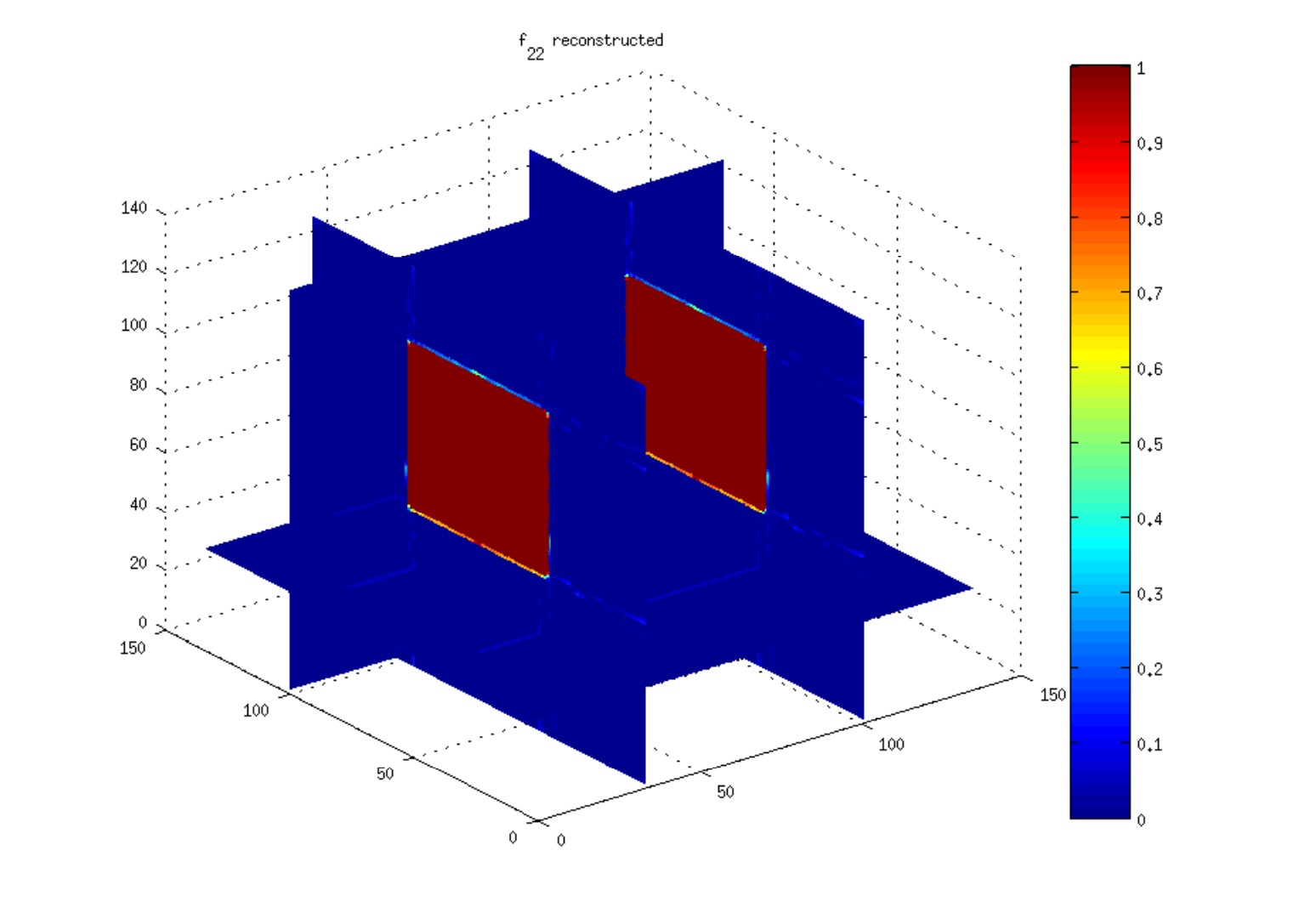}
        \caption{Reconstruction of $f_{22}$}
	\label{fig:rec_sharp_f22}
    \end{subfigure}
    \hspace{5mm}
    \begin{subfigure}[b]{0.45\textwidth}
        \centering
        \includegraphics[width=\textwidth]{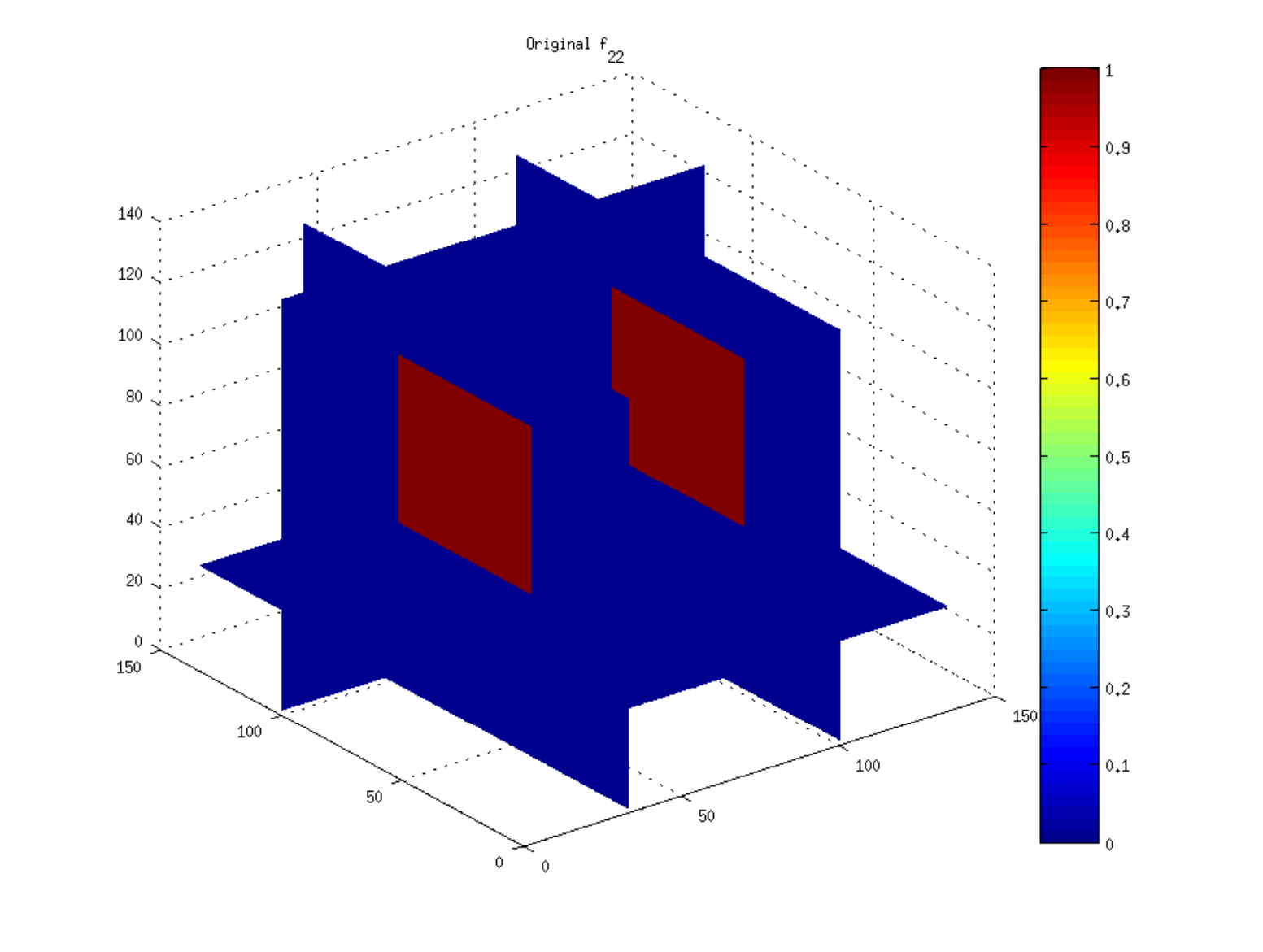}
        \caption{Original $f_{22}$}
	\label{fig:ori_sharp_f22}
    \end{subfigure}
    \vspace{5mm}
    \begin{subfigure}[b]{0.45\textwidth}
        \centering
        \includegraphics[width=\textwidth]{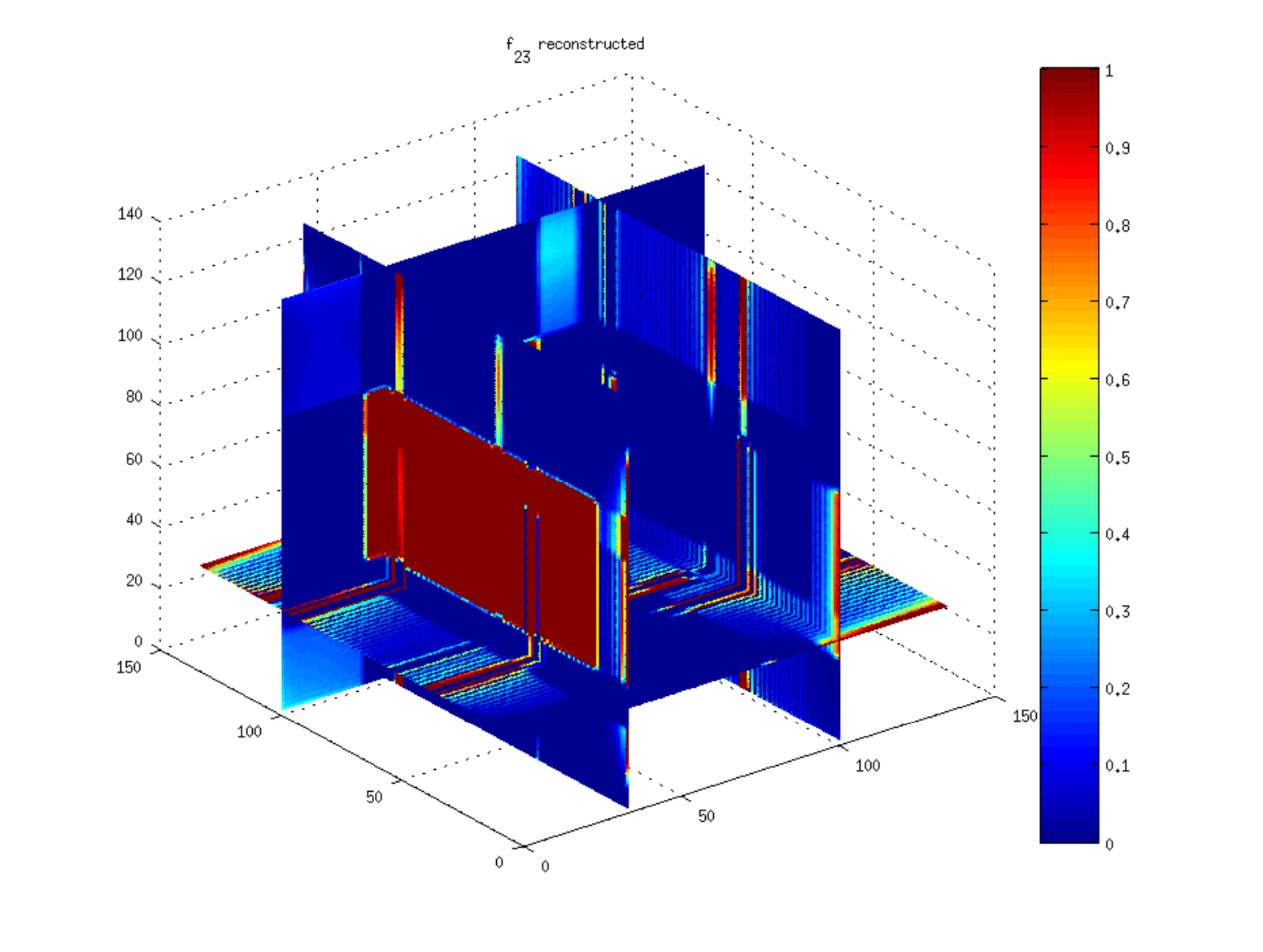}
        \caption{Reconstruction of $f_{23}$}
        \label{fig:rec_sharp_f23}
    \end{subfigure}
    \hspace{5mm}
    \begin{subfigure}[b]{0.45\textwidth}
        \centering
        \includegraphics[width=\textwidth]{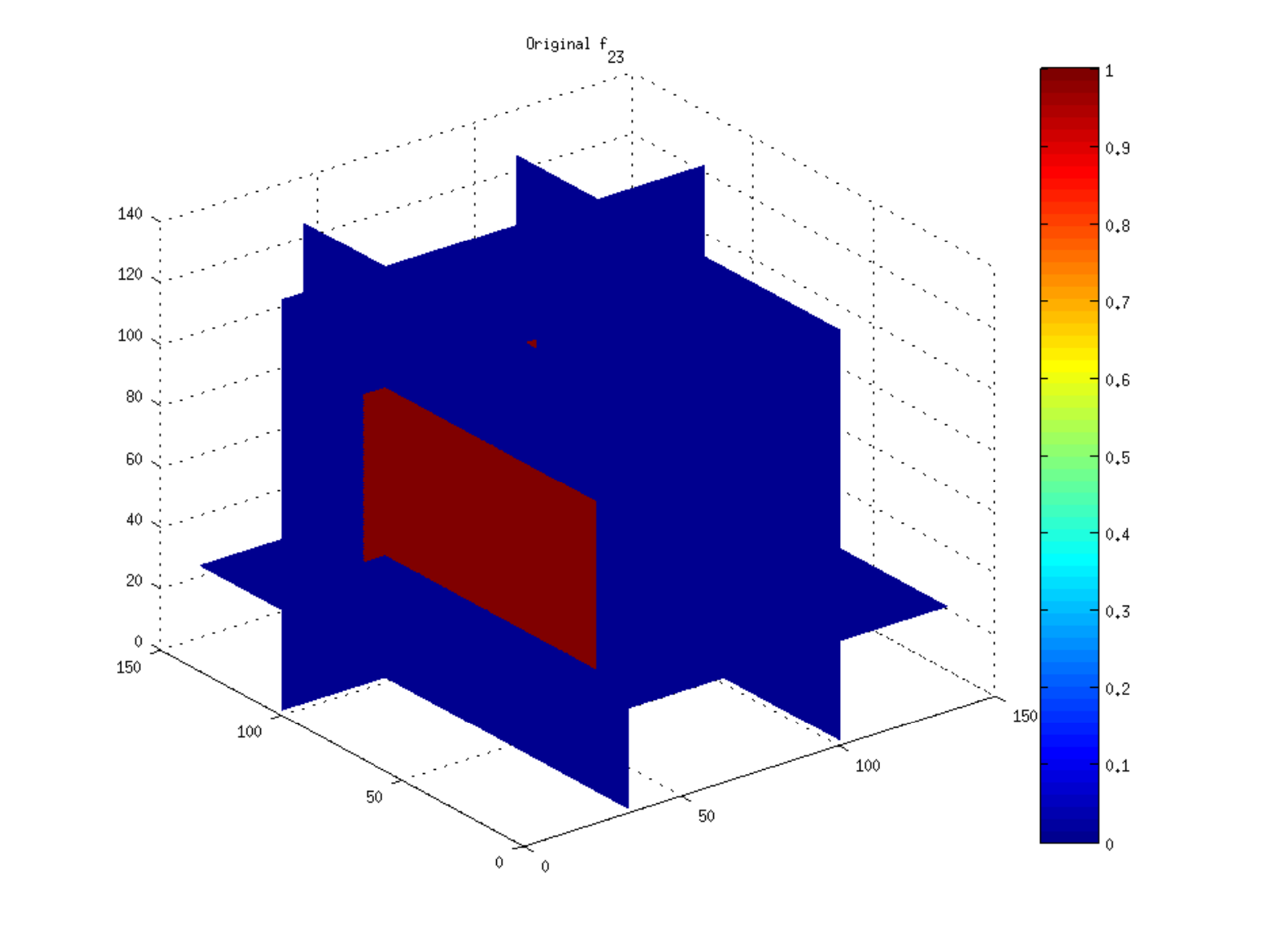}
        \caption{Original $f_{23}$}
        \label{fig:ori_sharp_f23}
    \end{subfigure}
    \vspace{5mm}
    \begin{subfigure}[b]{0.45\textwidth}
        \centering
        \includegraphics[width=\textwidth]{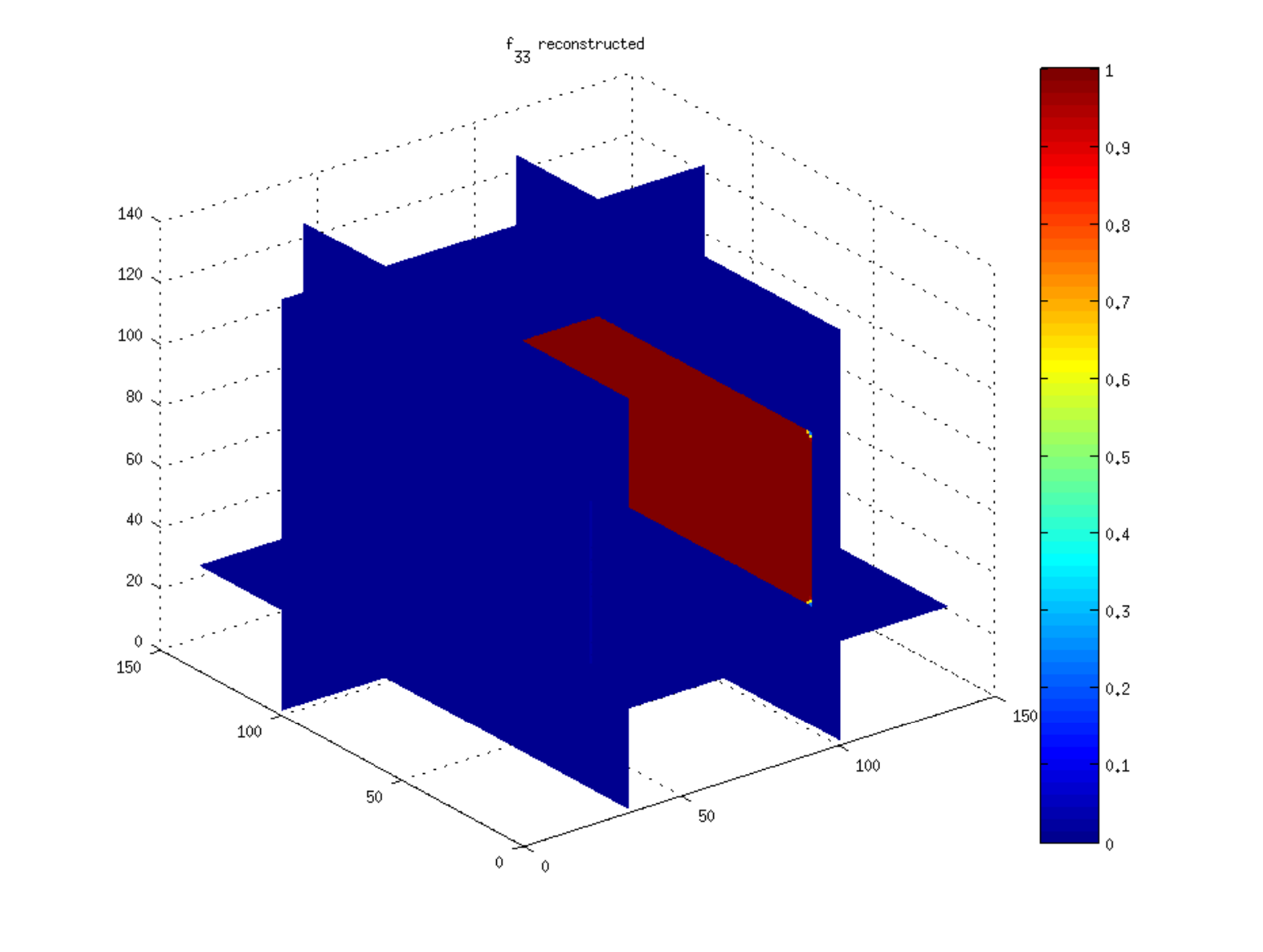}
        \caption{Reconstruction of $f_{33}$}
        \label{fig:rec_sharp_f33}
    \end{subfigure}
    \hspace{5mm}
    \begin{subfigure}[b]{0.45\textwidth}
        \centering
        \includegraphics[width=\textwidth]{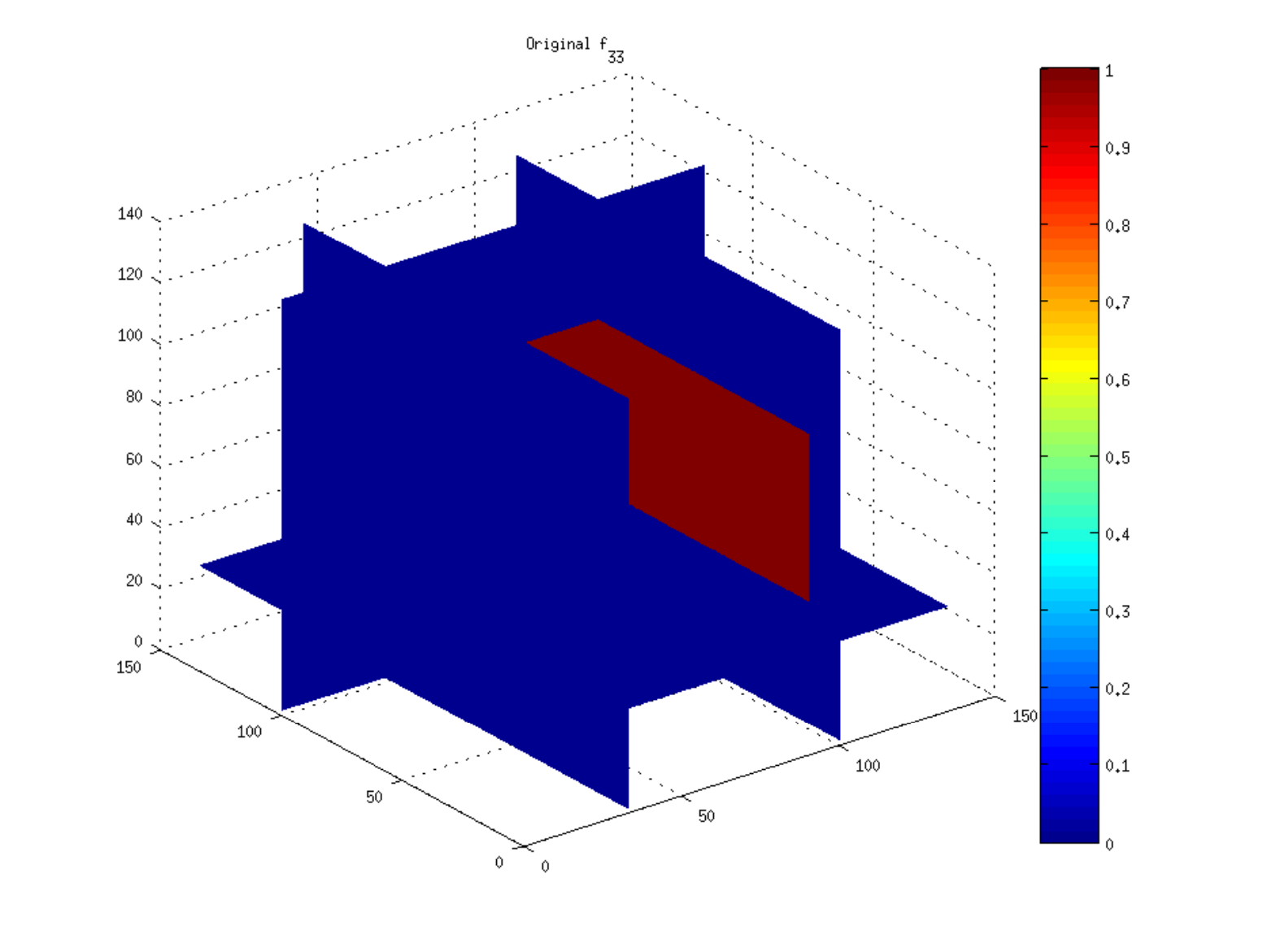}
        \caption{Original $f_{33}$}
        \label{fig:ori_sharp_f33}
    \end{subfigure}
    \caption{Reconstruction of sharp phantom, $f_{22}, f_{23}$ and $f_{33}$}
    \label{fig:sharp_f22_to_f33}
\end{figure}

\newpage
\section{Conclusions and further work}

Overall the procedure we describe would provide an efficient method for
x-ray diffraction tomography. It handles smooth tensor fields well but for discontinuous strain fields it would be sensible to
try different modified ramp filters or explicit regularization methods such as total variation.

If would be possible in the case of a tensor field that is the result of an infinitesimal strain to use only two axes of rotation, and recover the displacement field directly. However it might be better in practice to use the general procedure and then verify to what extent the compatibility condition holds on the reconstructed tensor.

We have pointed out that as there are two distinct methods of calculating the diagonal components this provides a consistency condition on the data. As the plane-by-plane data
is written in terms of scalar, vector and tensor longitudinal ray transforms, and the ranges of these operators can be determined in the plane case as a singular function expansion in a suitable Hilbert space.  In fan beam coordinates the singular value decomposition of the ray transform is given by \cite{Kazantsev}. See also \cite{Derevtsov2011} and  \cite{Derevtsov2014} for a parallel beam formulation. 

The  Helgason-Ludwig  range conditions for the scalar Radon transform in the plane \cite[Sec II.4]{Natterer} 
simply states that the $k$th moment of the data $\int\limits_{-\infty}^\infty s^k Xf(\xi,s\xi^\perp)\, \mathrm{d}s$, when it exists, is a polynomial of degree $\le k$ in $\xi$. For the LRT of a rank $m$ tensor the condition is simply degree $\le k+m$.  A deeper connection between these conditions and the singular function expansion is given by \cite{Monard}.

Such consistency conditions, characterizing the range of the forward operator are of great assistance in diagnosing errors and unaccounted for physical effects in experimental data.

Another avenue worth considering on the practical side  is to develop a reconstruction algorithm 
involving general (non-orthonormal) axes. In experiments it is often not feasible to rotate the specimen through $90^\circ$ and remain in the field of view of the measurement system. 

%It is important to consider  the speed of the actual reconstruction 
%when the data has been acquired, especially for large data-sets.  From \cite{DS}, it is apparent that inversion in Fourier space
%is a time consuming process.
%One possibility to speed up or compare the speed of reconstruction would be to switch to spatial domain and solve a 
%coupled system of partial differential equations represented as  finite difference approximations. 

Explicit reconstruction algorithms such as the one we have given are useful practically for data that is complete and uniformly sampled. For partial, sparse or irregularly sampled data
representing the forward problem simply as a sparse matrix and solving using iterative algorithms with explicit regularization is generally better, although typically requiring large amounts of memory and parallel processors.

\section*{Acknowledgements}
Parts of this work were supported by COST Action MP 1207, the Otto M{\o}nsteds Fond, the Royal Society Wolfson Research Merit Award, and EPSRC  Grants EP/M022498/1, % Tomographic Imaging CCPi
EP/M010619/1, %Next Generation Multi-Dimensional X-Ray Imaging
EP/I01912X/1, %MAPLE
P/K00428X/1  %MIRAN
\bibliography{TRT3}
\bibliographystyle{plain}

\end{document}